\documentclass{article}
\usepackage[english]{babel}
\usepackage[utf8]{inputenc} 
\usepackage[T1]{fontenc} 
\usepackage[letterpaper,top=2cm,bottom=2cm,left=3cm,right=3cm,marginparwidth=1.75cm]{geometry}
\usepackage{amsmath,amsfonts}
\usepackage{mathtools} 
\usepackage{dsfont}
\usepackage{stmaryrd} 
\usepackage{graphicx}
\usepackage{xcolor}
\colorlet{gray}{black!80!white}
\colorlet{blue}{blue!50!black}
\RequirePackage[colorlinks,citecolor=blue!60!black,urlcolor=blue!60!black,linktocpage, linkcolor=red!60!black]{hyperref} %
\usepackage{authblk}
\usepackage{easyReview}
\usepackage{caption}
\usepackage{subcaption}
\usepackage{amsthm} 
\theoremstyle{definition}
\newtheorem*{thm}{Theorem}
\newtheorem{theorem}{Theorem}[section]

\newtheorem{lemma}[theorem]{Lemma}

\theoremstyle{remark}
\newtheorem{remark}{Remark}[section]

\begin{document}
\title{A novel approach to the giant component fluctuations}
\author[1]{Josué Corujo}
\author[2]{Sophie Lemaire}
\author[3]{Vlada Limic}
\affil[1]{
Univ Paris Est Créteil, Univ Gustave Eiffel, CNRS, LAMA UMR 8050, F-94010 Créteil, France
}
\affil[2]{
Université Paris-Saclay, CNRS, Laboratoire de mathématiques d'Orsay, 91405, Orsay, France
}
\affil[3]{Institut de Recherche Mathématique Avancée, UMR 7501 Université de Strasbourg et CNRS, 
7 rue René-Descartes, 67000 Strasbourg, France}
{
	\makeatletter
	\renewcommand\AB@affilsepx{: \protect\Affilfont}
	\makeatother
	
	\affil[ ]{Email ids}
	
	\makeatletter
	\renewcommand\AB@affilsepx{, \protect\Affilfont}
	\makeatother
	
	\affil[1]{\href{mailto:josue.corujo-rodriguez@u-pec.fr}{josue.corujo-rodriguez@u-pec.fr}}
	\affil[2]{\href{mailto:sophie.lemaire@universite-paris-saclay.fr}{sophie.lemaire@universite-paris-saclay.fr}}
	\affil[3]{\href{mailto:vlada@math.unistra.fr}{vlada@math.unistra.fr}}
}

\maketitle

\begin{abstract}
We present a novel approach to study the evolution of the size (i.e.\ the number of vertices) of the giant component of a random graph process. It is based on the exploration algorithm called simultaneous breadth-first walks, introduced by Limic in 2019, that encodes the dynamics of the evolution of the sizes of the connected components of a large class of random graph processes. 
We limit our study to the variant of the Erd\H{o}s\,--\,R\'enyi graph process 
$(G_n(s))_{s\geq 0}$ with $n$ vertices where an edge connecting a pair of vertices appears at an exponential rate 1 waiting time, independently over pairs.  
We first use the properties of the simultaneous breadth-first walks to obtain an alternative and self-contained proof of the functional central limit theorem recently established by Enriquez, Faraud and Lemaire in the super-critical regime ($s=\frac{c}{n}$ and $c>1$). Next, to show the versatility of our approach, we prove a functional central limit theorem in the barely super-critical regime ($s=\frac{1+t\epsilon_n}{n}$ where $t>0$ and $(\epsilon_n)_n$ is a sequence of positive reals that converges to 0 such that $(n\epsilon_n^3)_n$ tends to $+\infty$). 
\end{abstract}

\smallskip


\smallskip

\section{Introduction and main results}

The Erd\H{o}s--R\'enyi random graph $\mathrm{ER}(n, p_n)$ has  $n$ vertices, and each pair of vertices is connected  by an edge with probability $p_n$,  independently over different pairs.
This model goes back to Gilbert \cite{1959Gilbert}, and a somewhat different version of the model was thoroughly analyzed by Erd\H{o}s and R\'enyi in their groundbreaking paper \cite{1959ErdosRenyi}.

Two vertices of $\mathrm{ER}(n, p_n)$ are \emph{related} if they are connected by a path in
$\mathrm{ER}(n, p_n)$ in the usual graph theoretical sense.
The above relation is clearly an
equivalence relation, and its  equivalence classes are called the \emph{connected components} of $\mathrm{ER}(n, p_n)$.
We refer to the number of vertices in a particular connected component $\mathcal{C}$
as the \emph{size of $\mathcal{C}$}.

\subsection{The phase transition and the super-critical regime}
\label{S:supercrit}
Let us take $p_n = c/n$.
When viewed on this scale, with varying $c\in [0,\infty)$, the model exhibits a phase transition at the critical parameter $c = 1$.
This phenomenon has been astonishing the probabilistic community ever since it was discovered by Erd\H{o}s and R\'enyi in their original paper \cite{1959ErdosRenyi}.
More precisely, if 
$c < 1$, then with overwhelming probability as $n\to \infty$,  all the connected components of the graph have size of order $O(\ln(n))$.
On the other hand, if $c = 1$ then for each fixed $k \in \mathbb{N}$, the $k$ largest connected components of $\mathrm{ER}(n, 1/n)$ all have sizes of order $\Theta(n^{2/3})$.
Finally, if $c > 1$, with overwhelming probability as $n\to \infty$ there is a single component of size
$\Theta(n)$, often called the \emph{giant component}, and all the other connected component sizes are $O(\ln(n))$.
This last phase is called the \emph{super-critical regime (phase)}. 
It is well-known that as $n\to \infty$ for a fixed $c>1$,  the largest component size is asymptotically equal to $\rho(c) n$, where $\rho(c)$ is the probability that a Bienaym\'e\,--\,Galton\,--\,Watson tree with a Poisson($c$) reproduction mechanism is infinite.
Moreover, $\rho(c)$ is the unique solution in $(0,1)$ of the equation
$\Phi^c (s) = 0$, where
\begin{equation}\label{def:Phi}
	\Phi^c : s \mapsto 1 - s - \mathrm{e}^{- c s}, \ s\geq 0.
\end{equation}

From now on we concentrate on the giant component.
Let us denote by $L_n^{\mathrm sc}(c)$ the size of the largest component in $\mathrm{ER}(n, c/n)$ with $c > 1$.
Stepanov \cite{Stepanov1970OnTP} proved  by a combinatorial method  
the following central limit theorem for $L_n^{\mathrm sc}(c)$:
\[
\sqrt{n} \left( \frac{L_n^{\mathrm sc}(c)}{n} - \rho(c) \right) \xrightarrow[n \to \infty]{\mathrm{d}} \mathcal{N}(0, \sigma^2(c)),
\]  where $ \sigma^2(c) = \frac{\rho(c)(1 - \rho(c))}{(1 - c (1 - \rho(c)))^2}$ and where $\xrightarrow{\mathrm{d}}$ denotes convergence in distribution. 
Other proofs of this result have been proposed 
by Pittel \cite{Pittel1990}, subsequently by Barraez, Boucheron and de la Vega~\cite{BBFdlV}  and Puhalskii~\cite{Puhalskii} via
exploration processes,
and more recently by R\'ath \cite{Rath2018} using moment generating functions.

\subsection{The barely super-critical regime}
\label{S:barely supercrit}
The barely super-critical happens ``in between the critical and the super-critical phases'', where there exists an ``early giant component'' with size which is already asymptotically larger than that of the other components.
In order to understand this regime, we need to refine the scale of the above phase-transition.

For a fixed $t>0$, let us consider the sequence of random graphs ${\mathrm{ER}(n, (1 + t \epsilon_n)/n)}$, where $(\epsilon_n)_n$ is a sequence of positive reals such that $\epsilon_n  \to 0$ and $n \epsilon_n^3 \to \infty$.
In this case, the size
$L^{\mathrm{bsc}}_n(t)$
of the largest component is
asymptotically equal to
${\rho(1 + t\epsilon_n) n = 2 t\epsilon_n n + O(\epsilon_n^2 n)}$.
The second largest connected component is of size $o(\epsilon_n n)$. These results were first proved by Bollob\'as \cite{1984Bollobas} with the stronger assumption $n \epsilon_n^3/\log(n)^{3/2} \to \infty$ and then by {\L}uczak \cite{1990Luczak} using enumerative methods. 
Nachmias and Peres presented in  \cite{2007NachmiasPeres} a simpler proof via an exploration process; using our notation their statement is
\begin{equation}\label{eq:result_Nachmias_Peres}
	\frac{L^{\mathrm{bsc}}_n(t)}{2 t\epsilon_n n} \xrightarrow[n \to \infty]{\mathbb{P}} 1 \ \text{ and } \ \frac{t^2 \epsilon_n^2}{2 \log(n t^3 \epsilon_n^3)} L^{\mathrm{bsc}; 2}_n(t) \xrightarrow[n \to \infty]{\mathbb{P}} 1,
\end{equation}
where $\xrightarrow{\mathrm{\mathbb{P}}}$ denotes convergence in probability.
Finally, Pittel and Wormald \cite{2005PittelWormald} proved that
\[
\sqrt{n \epsilon_n^3} \left( \frac{L_n^{\mathrm{bsc}}(t)}{n \epsilon_n} - \frac{\rho(1 + t \epsilon_n)}{\epsilon_n} \right) \xrightarrow[n \to \infty]{\mathrm{d}} \mathcal{N}\left(0, \frac{2}{t} \right),
\]
Bollob\'as and Riordan \cite{2012BollobasRiordan} proposed a simpler proof of the latter result, via an exploration process similar to the one of \cite{2007NachmiasPeres}.

\subsection{Motivation and main results}
All the above mentioned results focus on a fixed parameter ($c$ or $t$) within the super-critical or the barely super-critical phase.
However,
for each $n$, the size of the giant connected component can also be viewed as a process in $c$, via the natural (percolation-like) coupling of 
$(\mathrm{ER}(n, c/n), c>1)$.
Then for any $c_0>1$, and all large $n$, the process 
$
(L_n^{\mathrm{sc}}(c), \, c\geq c_0)
$
traces the size of the giant component.
Similarly, for any $t_0>0$, 
and all large $n$, the process 
$
(L_n^{\mathrm{bsc}}(t), \, t\geq t_0)
$
traces the size of the early giant component.

Enriquez, Faraud and Lemaire \cite{enriquez2023} recently initiated the study of giant component fluctuations in the context of varying $c$.
They proved the following functional limit theorem for the size of the largest component in the supercritical regime:
\begin{thm}[Theorem 1.1 \cite{enriquez2023}]
	Fix any $c_1>c_0 > 1$ and let $B$ denote the standard Brownian motion. 
	As $n\to \infty$,
	\[
	\left( \sqrt{n} \left( \frac{L_n^{\mathrm{sc}}(c)}{n} - \rho(c)\right),\, c \in [c_0, c_1]\right)
	\]
	converges weakly  to 
	\[
	\left( \frac{1 - \rho(c)}{1 - c (1 - \rho(c))} B \left( \frac{\rho(c)}{1 - \rho(c)} \right), \, c\in [c_0, c_1]\right),
	\] 
	in the usual Skorokhod topology. 
\end{thm}
The approach in \cite{enriquez2023} to prove this theorem was based on a detailed analysis of the infinitesimal increment of the size of the largest component, that led to a stochastic differential equation, which the authors of \cite{enriquez2023} explicitly solved.
In this article, we provide an alternative and shorter proof of this theorem using a novel method based  on a random  encoding for the sizes of the connected components of the random graph ER$(n, p=1-\mathrm{e}^{-t})$ as $t$ varies, called  \emph{the simultaneous breadth-first walks} by Limic in \cite{Limic2019}.   
Furthermore, we show the versatility of our technique by deriving a new result in the barely super-critical case. 

In order to state our results, let us describe  the coupling of the random graphs ER$(n,p)$, as $p$ varies, which 
is slightly different from that of \cite{enriquez2023}:
as in \cite{enriquez2023}, the initial state of the process is the completely isolated graph on $n$ vertices;
differently from \cite{enriquez2023}, we assume that an edge between the vertices $i$ and $j$ arrives at an exponential time, independently over pairs $(i,j)$, $i< j$. 
If the rate of the above exponential waiting times is $1$, we obtain the 
graph valued time-homogeneous Markov process (continuous-time Markov chain), denoted by $(\mathcal{G}^{(n)}(t))_{t \ge 0}$. 
Its marginal law at time $t\geq 0$ is the law of $\mathrm{ER}(n, p = 1 - \mathrm{e}^{-t})$.

Let us make an asymptotically negligible change in the notation used in Section \ref{S:supercrit}, so that now $L_n^{\mathrm{sc}}(c)$ and $L_n^{\mathrm{bsc}}(t)$ denote the size of the largest component in  $\mathcal{G}^{(n)}(c/n)$ and $\mathcal{G}^{(n)}((1 + t \epsilon_n)/n)$ respectively. 
Furthermore, let
$\Pi$ be a standard Brownian bridge on $[0,1]$, let $B$ denote the standard Brownian motion
and define a Gaussian process $(X(c), c >1)$  by
\begin{equation}\label{def:process_limit}
	X(c) = \frac{1}{1 - c (1 - \rho(c))} \Pi (\rho(c)).
\end{equation}
For the random graph process $(\mathcal{G}^{(n)}(t))_{t \ge 0}$, Theorem 1.1 \cite{enriquez2023} can be rewritten as follows: 
\begin{theorem}[Super-critical case - Theorem 1.1 \cite{enriquez2023}] \label{thm:enriquez}
	Let $L_n^{\mathrm{sc}}(c)$ denote the size of the largest component in  $\mathcal{G}^{(n)}(c/n)$. 
	Fix any $c_1>c_0 > 1$. 
	Then, as $n\to \infty$
	\[
	\left( \sqrt{n} \left( \frac{L_n^{\mathrm{sc}}(c)}{n} - \rho(c)\right),\, c \in [c_0, c_1]\right)
	\]
	converges weakly  to $(X(c), \, c\in [c_0, c_1])$ 
	in the usual Skorokhod topology. 
\end{theorem}

\begin{remark}
	The definition \eqref{def:process_limit} may seem quite different from 
	\[
	X'(c):=\frac{1 - \rho(c)}{1 - c (1 - \rho(c))} B \left( \frac{\rho(c)}{1 - \rho(c)} \right),
	\]
	which is the limit given in the statement of \cite[Theorem 1.1]{enriquez2023}, but it is easy to see that $X$ and $X'$ have identical laws.
	We prefer to state our result in terms of $\Pi$, which appears naturally in the description of the limit of the fluctuations of the simultaneous breadth-first walks. 
\end{remark}
The definition of simultaneous breadth-first walks and the context in which they appeared are presented in Section \ref{S:sBFW_and_MC}. 
Section \ref{S:pf_thm_1} is devoted to showing how to use simultaneous breadth-first walks to obtain a new proof of Theorem \ref{thm:enriquez}.  
In Section \ref{S:pf_thm_2} we prove the following result on 
the fluctuations of the size of the early giant component using an analogous approach.
\begin{theorem}[Barely super-critical case] \label{thm:barely-critical}
	Let $L_n^{\mathrm{bsc}}(t)$ denote the size of the largest component of $\mathcal{G}^{(n)}((1 + t \epsilon_n)/n)$. 
	Fix any $t_1 > t_0 > 0$. 
	If $(\epsilon_n)_n$ is a sequence of positive reals such that $\epsilon_n \to 0$ and $n \epsilon_n^3 \to \infty$ as $n\to \infty$, then
	\[
	\left( \sqrt{n \epsilon_n^3} \left( \frac{L_n^{\mathrm{bsc}}(t)}{n \epsilon_n} - \frac{\rho(1 + t \epsilon_n)}{\epsilon_n} \right), \, t\in [t_0, t_1] \right)
	\]
	converges weakly to the Gaussian process $\big( \frac{1}{t} B(2t), t \in [t_0, t_1] \big)$ in the usual Skorokhod topology.
\end{theorem}

\subsection{Conclusions and perspectives}
We believe that the use of simultaneous breadth-first walks has a considerable potential to yield additional functional limit theorems.
In fact, using similar techniques to those presented here, David Clancy, Jr.\ analyzes the fluctuations of inhomogeneous rank-one random graph processes \cite{Clancy_2025_rank1}, and the stochastic block model \cite{Clancy_2025_SBM}. 

A recent analysis by Bhamidi et al.~\cite{Bhamidi2024} of inhomogenous random graphs in the sub-critical and super-critical regimes relies on multi-type branching processes and stochastic analysis on Banach spaces. In the Erd\H{o}s-R\'enyi random graph case, their analysis provides an infinite dimensional functional central limit theorem for the sequence of densities of components of any given size. The limiting process is described as a solution of a stochastic differential equation. As a consequence, Bhamidi et al.~obtain in \cite[Theorem 3.11]{Bhamidi2024} the joint functional central limit theorem for the following three key quantities in the super-critical regime: the number of connected components, the size of the giant and the number of surplus edges in the giant component.  

We think that the joint fluctuations of at least the last two  quantities can also be obtained using the here presented approach. 
Indeed, it was recently shown that simultaneous breadth-first walks plus auxiliary randomness encode the dynamics of multiplicative coalescent with surplus edge counts (cf.\ \cite{Corujo_Limic_2023a}), thus enabling the study of the augmented multiplicative coalescents of Bhamidi et al.\ \cite{bhamidietal2} via natural and stronger tools. 
This insight suggests that the method we use to prove Theorem \ref{thm:enriquez} and Theorem \ref{thm:barely-critical} could be extended to encompass the fluctuations in the number of surplus edges.
In an ongoing project with Nathana\"el Enriquez and Gabriel Faraud we explore this research direction and one of our goals is to obtain a simple description of the joint scaling limit (as specified in \cite[Theorem 3.11]{Bhamidi2024}).

\section{The simultaneous breath-first walks and the multiplicative coalescent}
\label{S:sBFW_and_MC}

When we consider the evolution of connected component sizes of Erd\H{o}s-R\'enyi random graph evolving in continuous time (where each edge appears after an exponential time with rate $1$), we obtain the following dynamics:
\begin{quote}
	\emph{any two components of sizes  $x$ and $y$  merge at rate  $x \cdot y$ into a component of size $x + y$},
\end{quote}
which is called the \emph{multiplicative coalescent} dynamics since the seminal work of Aldous \cite{Aldous1997}, 

This motivates the following generalization of the random graph model.
Assume that the vertices have masses (in complete analogy with vertex size equal to $1$ in the previous section), and that the mass of vertex $i$ is $x_i > 0$.
Assume that the (undirected) edge connecting vertices $i$ and $j$ appears at an exponential time with rate $x_i \cdot x_j$, independently over $i\neq j$.
At each merger event, the mass of the new connected component is defined to be the sum of the masses,  prior to the merger, of the two components which have just coalesced.
Note that the evolution of the connected components masses in this model is again the multiplicative coalescent.
More precisely,
suppose that we are given an initial vector $x=(x_1, x_2, \dots, x_n)$ of vertex masses, such that $x_i \ge x_{i + 1}$.
For each $t$, let us denote by $X(t, x)$ 
the non-increasing list of connected component masses of the corresponding random graph at time $t$.
Then  $(X(t, x), t \ge 0)$ evolves according to the above (multiplicative coalescent) dynamics.

In  \cite{Aldous1997}, Aldous introduces a random function encoding based on the breadth-first algorithm used to explore  connected components of a graph in order to study $(X(t, x), t \ge 0)$ in the near-critical regime where all the masses of the components are of the same order. In particular, he showed that the weak limit at a fixed  $t$ of  $X(t+n^{1/3}, (n^{-2/3},\ldots,n^{-2/3}))$  as $n$ tends to $+\infty$ is the law of ordered excursion (above past minima) lengths of a Brownian motion with a parabolic drift. 
Since this seminal paper, different exploration processes for encoding the sizes of components of random graph were proposed and applied to extend this result to a model where the initial masses are not all of the same order (\cite{Aldous_Limic1998}) and after to prove a functional theorem as $t$ increases in different settings by Broutin and Marckert~\cite{BroutinMarckert},
Martin and R\'ath~\cite{Martin2017} and Limic~\cite{Limic2019} (see \cite{Martin2017} and \cite{Limic2019} for a comparison of these different exploration processes).   To the best of our knowledge, none of these exploration processes have been applied to prove a functional theorem outside the near-critical regime. 

The encoding proposed by Limic in \cite{Limic2019} is a modification of the breadth-first walk introduced in \cite{Aldous1997} in which the encoding functions are conveniently coupled so that their excursions above past infima encode the component sizes of $X(t,x)$ not only at a fixed coalescent time $t$, but also as $t$ varies.  It is defined with the help of a sequence of independent random variables $(\xi_i)_{i}$ such that $\xi_i$ is Exponential(rate $x_i$) distributed for each $1\leq i\leq n$, as follows:  
\[
Z^{x, t}: s \mapsto \sum_{i = 1}^n x_i \mathds{1}_{\{\xi_i/t\, \le\, s\}} - s, \text{ where } t>0. 
\]
The family of processes $\{Z^{x,t}, t\geq 0\}$ is called   the \emph{simultaneous breadth-first walks (sBFW)}. 

Recall that a segment $[A,B]$ is called an \emph{excursion} of
$Z^{x,t}$ \emph{(above past infima)} if, almost surely,
\begin{equation*}
	\label{def:excursion}
	Z^{x,t}(A-) = Z^{x,t}(B) = \inf_{s\leq A} Z^{x,t}(s), \ \ \text{ and } \ Z^{x,t}(s-)>Z^{x,t}(B), \ \forall s\in (A,B).
\end{equation*}
In the sequel, we simply write \emph{excursion} instead of excursion above past infima. 

We denote by $\Xi(t, x)$ the vector of lengths of the excursions of $Z^{x,t}$, where the excursion lengths are ordered non-increasingly and take  $\Xi(0, x) = x$. 
In \cite[Prop.~5]{Limic2019}, it is proved that 
the processes $(X(t, x), t \ge 0)$ and $(\Xi(t, x), t \ge 0)$ are equal in law.
It is therefore plausible that any result concerning the random graph connected component sizes (in fact, the multiplicative coalescent), can be derived using the two-parameter family of processes $\big( Z^{x,t}(s), s\geq 0, t > 0 \big)$.
In the next section, we show that in the context of results from \cite{enriquez2023}, this derivation is in fact shorter and simpler.

\section{The super-critical case}
\label{S:pf_thm_1}
We start by deriving the above-mentioned asymptotic for the size of the largest connected component in a supercritical regime, initially for a fixed $c>1$, and then uniformly on a segment of $(1,\infty)$.
\subsection{The concentration of the giant component size}
\label{S:concentration_sc}
\begin{figure}[htb]
	\centering \includegraphics[width=0.7\textwidth]{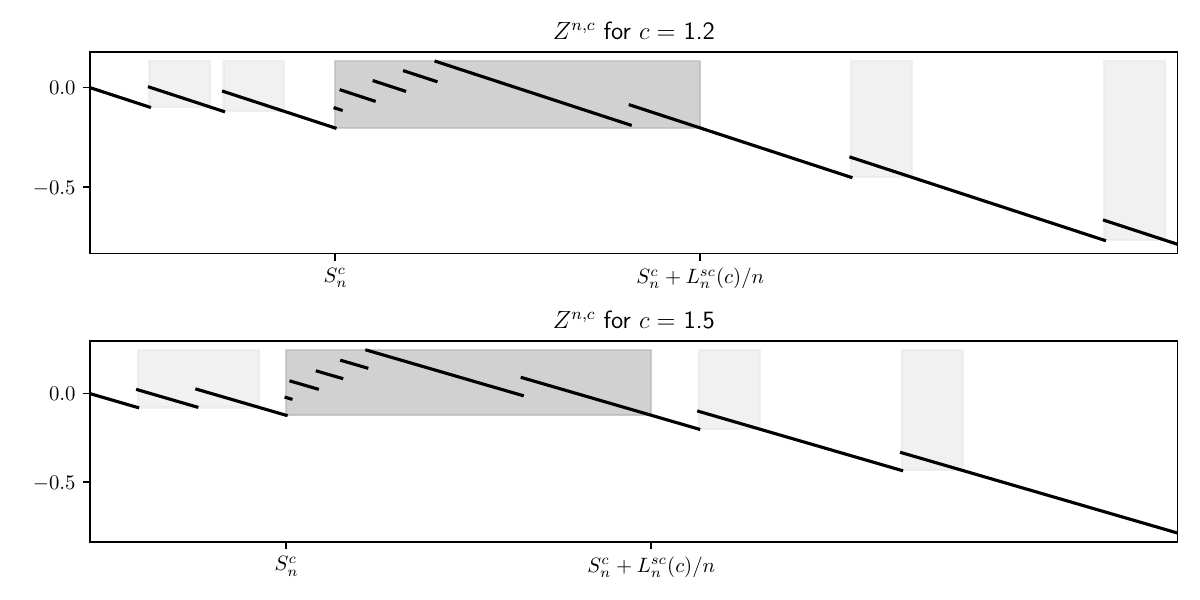}
	\caption{\label{fig_sBFWtraj} Evolution of a realization of the sBFW $s\mapsto Z^{n,c}(s)$ for $n=10$ as $c$ increases. 
		Excursions of $Z^{n,c}$ occur in the dark gray area for the longest excursion and in light gray areas for the other ones. The component sizes of the associated realization of the random graph  $\mathcal{G}^{(n)}(\frac{c}{n})$ are equal to $n$ times the interval lengths of the corresponding excursions of $Z^{n,c}$. In this example,  the associated graph $\mathcal{G}^{(n)}(\frac{c}{n})$ for $c=1.2$ has four isolated vertices and a component of size 6. As $c$ increases, the first two excursions merge;  the time at which they merge corresponds to the time at which a new edge is added in the random graph process $c\mapsto \mathcal{G}^{(n)}(\frac{c}{n})$ between the first two isolated vertices which are explored by $Z^{c}_{n}$. Between $c=1.2$ and $c=1.5$,  the largest excursion of $Z^{n,c}$ begins at  $S^{c}_n=\bar{\xi}_{(3)}/c$ where $\bar{\xi}_{(3)}$ denotes the third order statistic of $(\bar{\xi_i})_{1\leq i\leq n}$.}
\end{figure}

Consider a generalized random graph on $n$ vertices, where the  mass (size) of each vertex is $1/n$.
This implies that the mass of the giant component will be the first (and the only) coordinate (of order $1$) of $X(t_c^n, (1/n,\ldots,1/n))$, provided that $t_c^n$ is the multiplicative coalescent time which corresponds to $c/n$ in the random graph model $(\mathcal{G}^{(n)}(t))_{t \ge 0}$.
Note that the evolution of $X(t, (1/n,\ldots, 1/n))_{t\geq 0}$ is the same as the evolution of $(\mathcal{G}^{(n)}(s))_{s \ge 0}$ under the time change $s = t/n^2$ (the mass of each component in
$\mathcal{G}^{(n)}(t/n^2)$
is $n$ times larger than the mass of the corresponding component in
$X(t, (1/n,\ldots, 1/n))$). Therefore $t_c^n / n^2 = c/n$, or equivalently, $t_c^n = cn$.
Notice that if a random variable $\xi$ has a $\mathrm{Exponential}(\text{rate  } 1/n)$ distribution, then  $\bar{\xi} := \xi/n$ follows a $\mathrm{Exponential}( \text{rate  } 1)$ distribution.
Combining the above observations, we conclude that the relative size $L_n^{\mathrm{sc}}(c)/n$ of the largest component in the Erd\H{o}s\,--\,R\'enyi random graph equals in law to the length of the longest (largest) excursion  of
\[
Z^{n,c}(s) := \frac{1}{n} \sum_{i = 1}^n \mathds{1}_{\{ \bar{\xi}_i \le s c \}} - s,
\]
where $\bar{\xi}_i$ are i.i.d.\ with Exponential(rate  $1$) distribution.
See Figure \ref{fig_sBFWtraj} for a plot of an evolution of $Z^{n,c}$  as $c$ increases. 

Due to Glivenko\,--\,Cantelli theorem  we get that $(Z^{n,c}(s), s\ge 0)$ converges uniformly on $[0, \infty)$ a.s.\ (and therefore also in law in the Skorokhod $J_1$ topology) to $\Phi^c$ from \eqref{def:Phi}. 
We record this useful fact in symbols: 
\begin{equation}
	\label{equ:Z_unif_conv_Phi}
	\sup_{s \ge 0} | Z^{n, c}(s) - \Phi^c(s) | \xrightarrow[n \to \infty]{\text{a.s.}} 0.
\end{equation}

Our setting is similar to but much simpler than that of  \cite[Lemma 7]{Aldous1997} and \cite[Proposition 17]{Aldous_Limic1998}, since for any fixed $c>1$, $\Phi^c$ begins with an excursion of length $\rho(c)$ and is  strictly decreasing afterwards, with slope that we can control uniformly on $[c_0, \infty)$ for any $c_0>1$.
In particular, 
the uniform convergence \eqref{equ:Z_unif_conv_Phi}  guarantees that the longest excursion of 
$Z^{n,c}$ begins early (at the time $S_n^c$ converging to  $0$, as $n\to \infty$), that the length of this excursion  converges to $\rho(c)$ as $n\to \infty$, and that the length of the second longest excursion of $Z^{n,c}$ is negligible for all large $n$. 
In other words,
the starting point and the length of the longest excursion of $Z^{n, c}$ converge
to the starting point and the length of the longest (and unique) excursion of $\Phi^c$.
One can use real-analysis $\omega$-by-$\omega$ to verify the just stated asymptotics. 
More precisely, since the ``past infimum'' operator $\operatorname{I}$, defined by $\operatorname{I}(f): s \mapsto \inf_{u\leq s} f(u)$, is 1-Lipschitz continuous in the uniform topology, we get that
if for two functions $f$ and $f_{\text{app}}$ we have
$\sup_{s\geq 0} |f(s) - f_{\text{app}}(s)| \leq \delta$, then
it must be 
$\sup_{s\geq 0} | \operatorname{I}(f)(s) - \operatorname{I}(f_{\text{app}})(s)| \leq \delta$,
and therefore 
\begin{align}
	& 
	\forall \, s> 0 \text{ s.t.\ } f(s) - \operatorname{I}(f)(s) > 2 \delta \text{ it must be }f_{\text{app}}(s) - \operatorname{I}(f_{\text{app}})(s) >0, \label{equ:tall_excursion_f_implies_excursion_g}\\
	& 
	\forall s,a>0\text{ s.t.\ } f(s+a) < \operatorname{I}(f)(s)- 2\delta
	\text{ it must be }f_{\text{app}}(s+a)< \operatorname{I}(f_{\text{app}})(s). \label{equ:decrease_f_implies_no_excursion_g}   
\end{align}
From \eqref{equ:tall_excursion_f_implies_excursion_g} we conclude that whenever $f- \operatorname{I}(f)$  takes large enough (more than $2\delta$) values during an interval $[l,r]$, then  $f_{\text{app}}$ will necessarily have an excursion straddling $[l,r]$. 
Similarly, from \eqref{equ:decrease_f_implies_no_excursion_g}  
we learn that a significant decrease (of magnitude larger than $2\delta$) of $\operatorname{I}(f)$ over an interval $[s,s+a]$ prohibits existence of an excursion of $f_{\text{app}}$  straddling $[s,s+a]$.

We will soon address the vanishing of the second longest excursion length of $Z^{n,c}$ as $c$ varies.

As noticed by Enriquez et al.\ \cite[p.\ 2, eq.\ 4]{enriquez2023}, the quantity $\rho(c)$ defined on $[0,\infty)$, which satisfies $1 - \rho(c) = \mathrm{e}^{- c \rho(c)}$, is related to the the principal branch of the \emph{Lambert $W$ function}, denoted $W_0$, as follows
\[
\rho(c) = 1 + \frac{W_0(-c \, \mathrm{e}^{-c})}{c}.
\]
Moreover, $\rho$ is strictly increasing, continuous and converges to 1  as $c\to\infty$.

Since $L_n^{\mathrm{sc}}(c)/n$ is also a non-decreasing process in $c$, we get a stronger form of convergence via the following version of Dini's theorem: 
\begin{quote}
	\emph{If $(f_n)$ is a sequence of monotone functions (not necessarily continuous) converging pointwise to $f$ on $[a,\infty)$, and $f$ is continuous and bounded, then the convergence is uniform.}
\end{quote}
This yields
\begin{equation}
	\label{eq:uniform_as_conv}
	\sup_{c \ge c_0 }\left|\frac{L_n^{\mathrm{sc}}(c)}{n} - \rho(c) \right| \xrightarrow[n\to \infty]{\text{a.s.}}0,\quad \forall c_0>1.  
\end{equation}

Let us denote by $L_n^{\mathrm{sc};2}(c)$ the mass 
of the second largest  connected component of $\mathcal{G}^{(n)}(c/n)$.
Due to the sBFW encoding of our random graphs, we can and will also denote by $L_n^{\mathrm{sc};2}(c)/n$
the length of the second longest excursion of $Z^{n,c}$.
Note that $c\mapsto L_n^{\mathrm{sc};2}(c)$
is not a monotone map, almost surely.
We will need  some control over $L_n^{\mathrm{sc};2}(c)$ as $c$ varies.  
Let $\operatorname{id}$ be the identity map on $\mathbb{R}$, and note
that $Z^{n,c} + \operatorname{id}$ and $\Phi^c + \operatorname{id}$ scale compatibly with respect to $c$, so that for every $c>c_0>0$ we have
\[
Z^{n,c}(s) - \Phi^c(s)  =  Z^{n,c_0} \left( \frac{c}{c_0} s \right) - \Phi^{c_0} \left( \frac{c}{c_0} s \right), \ \ s\geq 0.
\]
Therefore 
\[
\sup_{s \ge 0} | Z^{n, c_0}(s) - \Phi^{c_0}(s) |
=
\sup_{c \ge c_0}\, 
\sup_{s \ge 0} | Z^{n, c}(s) - \Phi^c(s) |.
\]
Focusing on the supercritical phase, we fix $c_0>1$. We get from \eqref{equ:Z_unif_conv_Phi} that  
\begin{equation}\label{equ:Z_unif_conv_Phi_uniform_c}
	\sup_{c \ge c_0}\, 
	\sup_{s \ge 0} | Z^{n, c}(s) - \Phi^c(s) |\
	\xrightarrow[n \to \infty]{\text{a.s.}} 0,
\end{equation}
Note that $(\Phi^c)'(\rho(c)) = \pi(c) - 1$, where  $\pi(c) := c (1 - \rho(c))$  is the conjugate to $c$ in the sense that 
$\pi(c) \mathrm{e}^{-\pi(c)} = c \, \mathrm{e}^{-c}$ and  $\pi(c) < 1$  when $c > 1$ (cf.\ e.g.~\cite[Thm.\ 3.15]{Remco2017}).
So we get a uniform bound
\[
\sup_{c \ge c_0} \sup_{s \ge \rho(c)} (\Phi^c)'(s) = (\Phi^{c_0})'(\rho(c_0)) = \pi(c_0) - 1 < 0,\quad \forall c_0>1.
\]
Joint with a simultaneous application of \eqref{equ:tall_excursion_f_implies_excursion_g}--\eqref{equ:decrease_f_implies_no_excursion_g} 
with $(f,f_{\text{app}})=(\Phi^c,Z^{n,c})$ for all $c\in [c_0,\infty)$ 
this implies that the length of the second longest excursion of $Z^{n,c}$ is negligible uniformly in $c\geq c_0$, as $n\to \infty$.
To summarize, we have argued that
\begin{equation}
	\label{equ:second_largest_negligible}
	\sup_{c \geq c_0} \left| \frac{L_n^{\mathrm{sc};2}(c)}{n} \right| \xrightarrow[n\to \infty]{\text{a.s.}} 0.
\end{equation}

Furthermore,  the almost sure convergence \eqref{equ:Z_unif_conv_Phi}, \eqref{eq:uniform_as_conv} and  \eqref{equ:second_largest_negligible} imply that with probability $1$, for $n$ sufficiently large and uniformly in $c \ge c_0$, the longest excursion can be identified (as the one straddling the point $\rho(c_0)/2$).
Let us denote by $S_n^c$ the starting point of the longest excursion of $Z^{n,c}$.
Note that $S_n^c = \bar{\xi}_{(I_n^c)}/c$ for some random index $I_n^c$.
Since the longest excursion traverses $\rho(c_0)/2$ for all $n$ sufficiently large and all $c\geq c_0$, we conclude that 
$I^c_n$ is non-increasing in $c$ and therefore that
$(S_n^c, c \ge c_0)$ is (strictly) decreasing in $c$ for all large $n$ (see Figure \ref{fig_sBFWtraj} for an example of how $S^{c}_n$ changes as $c$ increases). 
Thus,
\begin{equation}
	\label{eq:uniform_as_S}
	\sup_{c\in [c_0, \infty)} S_n^c \xrightarrow[n\to \infty]{\text{a.s.}} 0.
\end{equation}
The control of the fluctuations in convergence \eqref{eq:uniform_as_conv} is a more interesting and harder problem.
As we already mentioned, it was recently solved in \cite{enriquez2023}, but one aim of this article is to exhibit an alternative and shorter proof.

\subsection{The fluctuations of the giant component size}
\label{S:proof_thm_enriquez}
Our goal is to prove that 
\[
\left( \sqrt{n} \left( \frac{L_n^{\mathrm{sc}}(c)}{n} - \rho(c) \right),\; c\in[c_0,c_1]\right),
\]
converges weakly by using the fact that $L_n^{\mathrm{sc}}(c)/n$ is the length of the longest excursion of $Z^{n, c}$. 
Let us note that while this object is well defined, there is no simple algebraic expression for $L_n^{\mathrm{sc}}(c)/n$ in terms of the input data $(\bar{\xi}_i)_i$. 
The idea is to analyze the fluctuations of  $\Phi^{(c)}(S^{c}_n+\frac{L_n^{\mathrm{sc}}(c)}{n})$ and apply the mean value theorem to $\Phi^c$ between $S_n^c+\frac{L_n^{\mathrm{sc}}(c)}{n}$, which is the endpoint of the longest excursion of $Z^{c}_n$ and $\rho(c)$, which is the endpoint of the unique excursion of $\Phi^c$. 
Let us first have a look at the behavior of the centered breadth-first walk
\[
\sqrt{n} \left( Z^{n, c}\big( s \big) - \Phi^c(s) \right) = \sqrt{n} \left( \frac{1}{n} \sum_{i = 1}^n \mathds{1}_{\{ \overline{\xi}_i \le c s \}} 
-  (1 - \mathrm{e}^{-c s}) \right)
\]
for a fixed $c>1$. 

Define $F^c(s) = 1 - \mathrm{e}^{- c s}$ to be the cumulative distribution function of $\overline{\xi}_i/c$.
Applying the version of Donsker's theorem for empirical processes (cf.\ \cite[Theorem ~1.1.1]{Dudley1999}) 
we get that 
\begin{lemma}[Donsker's theorem for $Z^{n, c}$]\label{lemma:Donsker}
	The process 
	\[
	\left( \sqrt{n} \left( Z^{n, c}\big( s \big) - \Phi^c(s) \right), s \ge 0 \right)
	\]
	converges weakly in $\mathcal{C}([0, \infty))$ towards $(\Pi \big( F^c(s) \big), s \ge 0)$, where $(\Pi(s), s \in [0,1])$ is a standard Brownian bridge.
\end{lemma}
Let us fix some $c_1>1$.
Due to Skorokhod's representation theorem, there is a family of random functions $(\hat{Z}^{n,c_1})_{n \ge 1}$ and a Brownian bridge $\hat{\Pi}$, defined both on some probability space, such that for every $n \ge 1$, $\hat{Z}^{n,c_1}$ is equal in law to ${Z}^{n,c_1}$, and for every $s_1 > 0$ we have
\begin{equation}
	\label{eq:kolmogorov}
	\sup_{s \in [0, s_1]} \left| \sqrt{n} \left[  \hat{Z}^{n,c_1} \left( s \right) - \Phi^{c_1} \left( s \right) \right] - \hat{\Pi} \big( F^{c_1}(s) \big) \right| \xrightarrow[n \to \infty]{\mathrm{a.s.}} 0.
\end{equation}
Note that, for each $n$, the equality in law between $\hat{Z}^{n,c_1}$ and $Z^{n,c_1}$ implies that
$$
\hat{Z}^{n,c_1}(s) = \frac{1}{n} \sum_{i = 1}^n \mathds{1}_{\{ \overline{\xi}^{(n)}_i \le c_1 s \}}-s,$$ 
where the unordered jump times  ${(\bar{\xi}_i^{(n)}/c_1)_{1\leq i\leq n}}$ are a family of i.i.d.~exponential rate $c_1$ random variables.
There is no a priori link between
${(\bar{\xi}_i^{(n_1)})_{1\leq i\leq n_1}}$
and
${(\bar{\xi}_i^{(n_2)})_{1\leq i\leq n_2}}$
if $n_1\neq n_2$.

For each $c > 0$, let us define
\begin{equation}
	\hat{Z}^{n,c}(s):=\frac{1}{n} \sum_{i = 1}^n \mathds{1}_{\{ \overline{\xi}^{(n)}_i \le c s \}}-s 
	= (\hat{Z}^{n,c_1} + \operatorname{id}) \left(\frac{sc}{c_1}\right) - s. 
\end{equation}
Notice that for every $n \ge 1$, we have 
\begin{equation}
	\label{eq:hatZ_ed_Z}
	(\hat{Z}^{n,c}(s), s\geq 0,\ c \ge 0) \stackrel{d}{=} ({Z}^{n,c}(s), s\geq 0,\ c \ge 0).
\end{equation}
Let us also define
\begin{equation}\label{eq:def_Gamma}
	\Gamma_n^c(s) := \sqrt{n} \left[  \hat{Z}^{n,c} \left( s \right) - \Phi^c \left( s \right) \right] - \hat{\Pi} \big( F^c(s) \big).
\end{equation}

In order to extend \eqref{eq:kolmogorov} to the dynamic setting,
we will use again the scaling compatibly with respect to $c$.
Indeed, since
$\hat{Z}^{n,c} - \Phi^c$ and $F^c$ scale compatibly, we have that for every $c_0\in (0,c_1)$ and all $c \in [c_0, c_1]$ and $s_1>0$
\[
\sup_{s \in [0, s_1]} | \Gamma^{c}_n(s) | = \sup_{s \in [0 , s_1 \frac{c}{c_1}]} | \Gamma^{c_1}_n(s) |.
\]
Thus, on the above probability space the convergence \eqref{eq:kolmogorov} is uniform: for every $c_0 \in  (1, c_1)$ and any $s_1>0$, we have
\begin{equation}
	\label{equ:unif_Psi_conv}
	\sup_{c \in [c_0, c_1] }  \sup_{s \in [0, s_1]} | \Gamma^c_n(s) | = \sup_{s \in [0, s_1]} | \Gamma^{c_1}_n(s) | \xrightarrow[n \to \infty]{\text{a.s.}} 0. 
\end{equation}
Let $\hat{S}^c_n$ and $\hat{L}_n^{\mathrm{sc}}(c)$ be the starting point and the length of the longest excursion of $\hat{Z}^{n,c}$, respectively.
As a consequence of \eqref{eq:hatZ_ed_Z}, 
\begin{equation}
	\label{eq:hatSL_eq_SL}
	(\hat{S}^c_n)_{c \in [c_0, c_1]}\stackrel{d}{=}(S^c_n)_{c \in [c_0, c_1]}\text{ and }(\hat{L}^{\mathrm{sc}}_n(c))_{c \in [c_0, c_1]}\stackrel{d}{=}({L}^{\mathrm{sc}}_n(c))_{c \in [c_0, c_1]}.
\end{equation}
Since $\Phi^c(\rho(c)) = 0$, the mean value theorem applied to $\Phi^c$ gives 
\begin{equation}\label{eq:Taylor_expansion_Phi}
	\Phi^c \left(\hat{S}_n^c + \frac{\hat{L}_n^{\mathrm{sc}}(c)}{n} \right)  = (\Phi^c)'(\alpha_n^c) \left( \hat{S}_n^c + \frac{\hat{L}_n^{\mathrm{sc}}(c)}{n} - \rho(c) \right)
\end{equation}
for some  $\alpha_n^c \in \big[ \rho(c) \wedge (\hat{S}_n^c + \frac{\hat{L}_n^{\mathrm{sc}}(c)}{n}), \rho(c) \vee (\hat{S}_n^c + \frac{\hat{L}_n^{\mathrm{sc}}(c)}{n}) \big]$.
Furthermore, 
\begin{align}
	\sup_{c\in [c_0,c_1]}|(\Phi^c)'(\alpha_n^c)-(\Phi^c)'(\rho(c))|
	&= 
	\sup_{c\in [c_0,c_1]} c |\mathrm{e}^{- c\alpha_n^c }- \mathrm{e}^{-c \rho(c)}| \nonumber \\ 
	&\le \sup_{c\in [c_0,c_1]} c^2 | \alpha_n^c - \rho(c) |
	\xrightarrow[n \to \infty]{\mathbb{P}} 0, \label{conv_phic_prime}
\end{align}
where the last convergence is a consequence of  \eqref{eq:uniform_as_conv},  \eqref{eq:uniform_as_S} and  \eqref{eq:hatSL_eq_SL}.

Moreover, for every $c > 1$ we have
\begin{equation}\label{eq:deriv_rho}
	(\Phi^{c})'(\rho(c)) = c(1-\rho(c)) -1 <0.
\end{equation}

Let us define $\gamma_n^c := \hat{S}_n^c + \frac{\hat{L}_n^{\mathrm{sc}}(c)}{n}$
and recall the definition of $\Gamma_n^c$ in \eqref{eq:def_Gamma}.
We can now rewrite the sum of the term on the left-hand side of \eqref{eq:Taylor_expansion_Phi} (multiplied  by $\sqrt{n}$) and $ \hat{\Pi} \big(F^c(\rho(c)) \big)$ as
\[
\sqrt{n} \Phi^c \left( \gamma_n^c \right) + \hat{\Pi} \big( F^c(\rho(c)) \big)
\]
which equals
\begin{equation}
	\label{eq:rhs_Taylor_expansion_Phi}
	- \Gamma_n^c(\gamma_n^c) 
	+ \sqrt{n} \hat{Z}^{n,c} \left( \gamma_n^c \right)  - \hat{\Pi} \big( F^c(\gamma_n^c) \big) + \hat{\Pi} \big( F^c(\rho(c)) \big).
\end{equation}
Due to \eqref{equ:unif_Psi_conv}, the first term on the right-hand side of the previous display converges to zero a.s.\ uniformly on compacts as $n \to \infty$. 
Furthermore, due to the uniform continuity of the Brownian bridge, $\hat{\Pi} \big( F^c(\rho(c)) \big) - \hat{\Pi} \big( F^c(\gamma_n^c) \big)$ also converges to zero uniformly on compacts.
The remaining term 
in \eqref{eq:rhs_Taylor_expansion_Phi} can be controlled using the first part of the following lemma.
\begin{lemma}\label{lemma:control}
	For every $1 < c_0 < c_ 1 < \infty$, the non-negative processes 
	\begin{equation*}
		\Bigg( -\sqrt{n} \hat{Z}^{n,c} \Big( \hat{S}_n^c + \frac{\hat{L}_n^{\mathrm{sc}}(c)}{n} \Big), c \in [c_0,c_1]  \Bigg) \text{ and }
		\Bigg( \sqrt{n} \hat{S}_n^c, c \in[c_0,c_ 1] \Bigg)
		\label{equ:lemma:control}
	\end{equation*}
	converge uniformly in probability to zero, as $n \to \infty$. 
\end{lemma}
Before proving Lemma \ref{lemma:control}, let us conclude the proof of Theorem \ref{thm:enriquez}.
Combining the observations made above with Lemma \ref{lemma:control}, we get that
\[
\sqrt{n} \Phi^c \left( \gamma_n^c \right) + \hat{\Pi} \big( F^c(\rho(c)) \big)
\]
converges uniformly on compacts in probability to zero, as $n \to \infty$.
Hence, recalling \eqref{eq:Taylor_expansion_Phi} and Lemma \ref{lemma:control}
we get that
\[
\big( \Phi^c \big)'(\alpha_n^c) \cdot \sqrt{n} \left( \frac{\hat{L}_n^{\mathrm{sc}}(c)}{n} - \rho(c) \right) + \hat{\Pi} \big( F^c(\rho(c)) \big)
\]
converges uniformly on compacts in probability  to zero.
Now we apply \eqref{conv_phic_prime} and \eqref{eq:deriv_rho} to obtain that
$$\left( \displaystyle{ \sqrt{n} \left( \frac{\hat{L}_n^{\mathrm{sc}}(c)}{n} - \rho(c) \right)}, c \in [c_0, c_1] \right)$$
converges uniformly in probability to
$$\left( \displaystyle{\frac{1}{1 - c (1 - \rho(c))} \hat{\Pi} \big( F^c(\rho(c)) \big)}, c \in [c_0, c_1] \right),$$
as $n$ tends to  $+\infty$.
Note that $F^c (\rho(c) ) = \rho(c)$.
This concludes the proof of Theorem \ref{thm:enriquez}, since  $(\hat{L}^{\mathrm{sc}}_n(c))_{c \in[c_0,c_1]}$ has the same distribution as $(L^{\mathrm{sc}}(c))_{c\in[c_0,c_1]}$ by Proposition 5 in \cite{Limic2019}.

\begin{proof}[Proof of Lemma \ref{lemma:control}]
	Recall that $Z^{n, c}$ is the ``original'' sBFW and $\hat{Z}^{n, c}$ is the sBFW in the space of the Skorokhod's representation.
	Then, for every $\epsilon > 0$ and each $n \geq 1$ we have
	$$
	\mathbb{P} \left( 
	\sup_{s \ge 0} | \hat{Z}^{n, c}(s) - \Phi^c(s) |\ \le \epsilon  \right) = \mathbb{P} \left(
	\sup_{s \ge 0} | {Z}^{n, c}(s) - \Phi^c(s) |\ \le \epsilon  \right).
	$$
	Furthermore, because of the scaling compatibility property we have that
	$$
	| {Z}^{n, c_0}(s) - \Phi^{c_0}(s) | = \left| {Z}^{n, c}\left(\frac{s c_0}{c}\right) - \Phi^{c}\left(\frac{s c_0}{c} \right) \right|.
	$$
	This implies that for each $n\geq 1$
	$$
	\left\{ 
	\sup_{s \ge 0} | Z^{n, c_0}(s) - \Phi^{c_0}(s) |\ \le \epsilon \right\} =  \left\{  \sup_{c \ge c_0}\, 
	\sup_{s \ge 0} | Z^{n, c}(s) - \Phi^c(s) |\ \le \epsilon \right\}
	$$
	and the same identities are true for $\hat{Z}^{n,c}$, almost surely.
	In the end, we obtain for each $n\geq 1$
	\[
	\mathbb{P} \left( \sup_{c \ge c_0}  
	\sup_{s \ge 0} | \hat{Z}^{n, c}(s) - \Phi^c(s) |\ \le \epsilon \right) = \mathbb{P} \left( \sup_{c \ge c_0}  
	\sup_{s \ge 0} | {Z}^{n, c}(s) - \Phi^c(s) |\ \le \epsilon  \right).
	\]
	The uniform convergence \eqref{equ:Z_unif_conv_Phi_uniform_c} together with the just stated identity gives 
	\[
	\lim_{n} \mathbb{P} \left( \sup_{c \ge c_0}   \sup_{s \ge 0} | \hat{Z}^{n, c}(s) - \Phi^c(s) |\ \le \epsilon  \right) = 1.
	\]
	Let us take $\epsilon=\epsilon(c_0)>0$ sufficiently small so that on
	\begin{equation}\label{eq:event_all_c}
		\left\{  \sup_{c \ge c_0}\, 
		\sup_{s \ge 0} | \hat{Z}^{n, c}(s) - \Phi^c(s) |\ \le \epsilon \right\}
	\end{equation}
	for all $c \geq c_0$, the longest excursion of $\hat{Z}^{n, c}$ is its only excursion longer than $\rho(c_0)/2$, and such that it straddles $[\rho(c_0)/4, 3 \rho(c_0)/4]$.  
	Therefore, 
	on the same event,
	the longest excursion left endpoint $\hat{S}_n^c$ is well defined for all $c\geq c_0$, and it is a decreasing function of $c$.
	Recall that the probability of this event converges to $1$ as $n\to \infty$. 
	Elementary properties of $\hat{Z}^{n,c}$ imply 
	$$
	\hat{S}_n^c \ge -\hat{Z}^{n, c}(\hat{S}_n^c-) = - \hat{Z}^{n, c} ( \hat{S}_n^c + \hat{L}_n^{\mathrm{sc}}(c)/n) \ge 0,
	$$
	so it suffices to show
	that $(\sqrt{n} \hat{S}_n^c)_n$ converges to $0$ in probability with $c=c_0$.
	
	Let us write $c$ for $c_0$ to reduce notational clutter.
	Recall that the left limit of $\hat{Z}^{n,c}$ is strictly negative at the beginning of each excursion (see Figure \ref{fig_sBFWtraj}).
	The idea is to prove that for any $\delta > 0$, with overwhelming probability as $n  \to \infty$, 
	the process $\hat{Z}^{n,c}$ takes non-negative values on $ [{\delta}/{\sqrt{n}},  \alpha_c]$, for some deterministic $\alpha_c>0$.
	This would imply that for $n$ large enough, with overwhelming probability, the interval $[{\delta}/{\sqrt{n}},  \alpha_c]$ must be contained in some excursion of $\hat{Z}^{n,c}$.
	As argued before, the almost sure convergence \eqref{equ:Z_unif_conv_Phi_uniform_c} implies that with high probability (w.h.p.) the excursion covering $[{\delta}/{\sqrt{n}},  \alpha_c]$ is the largest one.
	Thus, the event \eqref{eq:event_all_c} intersected with $\{ \hat{S}_n^c \le {\delta}/{\sqrt{n}} \}$ would have asymptotic probability one.
	This would prove that $\sqrt{n} \hat{S}_n^c$ converges to zero in probability.
	
	Let us now prove that $\hat{Z}^{n,c}$ is indeed positive on an interval of the form $[{\delta}/{\sqrt{n}},  \alpha_c]$, for some $\alpha_c > 0$.
	Lemma \ref{lemma:Donsker} yields the following asymptotics on compacts
	\[
	\hat{Z}^{n, c}\big( s \big) =  \Phi^c(s) + \frac{1}{\sqrt{n}} \hat{\Pi}(F^c(s)) + o_{\mathbb{P}}\left( \frac{1}{\sqrt{n}} \right).
	\]
	Let us choose $\alpha_c$ as the argmax of $\Phi^c$ in $[0,1]$, which is actually the point of its global maximum in $\mathbb{R}_+$, due to strict concavity.
	We conveniently split the above interval into two pieces: $[{\delta}/{\sqrt{n}}, {1}/{\sqrt[4]{n}}]$ and $[{1}/{\sqrt[4]{n}}, \alpha_c]$, where $\sqrt[4]{n}= \sqrt{\sqrt{n}}$.
	On the one hand, for all $s\in [{\delta}/{\sqrt{n}}, {1}/{\sqrt[4]{n}}]$ we have
	\[
	\sqrt{n} \hat{Z}^{n,c} (s) \ge \sqrt{n} \Phi^c \left( \frac{\delta}{\sqrt{n}} \right) + o_{\mathbb{P}} \left( 1 \right) \xrightarrow[n \to \infty]{\mathbb{P}} \delta \cdot (\Phi^c)'(0) = \delta (c - 1) > 0,
	\]
	{where we also use the continuity of the Brownian bridge at $0$.}
	On the other hand, for all $s \in [{1}/{\sqrt[4]{n}}, \alpha_c]$ we have
	\[
	\sqrt{n} \hat{Z}^{n, c}(s) \ge \sqrt{n} \Phi^c \left( \frac{1}{\sqrt[4]{n}} \right) + O_{\mathbb{P}} \left( 1 \right) = \Omega_{\mathbb{P}}(\sqrt[4]{n}),
	\]
	where we again use the differentiablity of $\Phi^c$ at $0$, as well as the boundedness of  Brownian bridge on compacts.
	It is easy to see that the last two estimates imply 
	\[
	\inf_{s \in [\frac{\delta}{ \sqrt{n}}, \alpha_c]} \hat{Z}^{n, c}(s) > 0 \text{ w.h.p}.
	\]
	As already argued, this concludes the proof of the lemma.
\end{proof}

\section{The barely super-critical case}
\label{S:pf_thm_2}

This section is devoted to the study of the barely super-critical regime: we consider the random graph process $(\mathcal{G}^{(n)}((1+t\epsilon_n)/n))_{t>0}$ where $t>0$ and $(\epsilon_n)_n$ is a sequence of positive reals such that $\epsilon_n\rightarrow 0$ and $n\epsilon^{3}_n\rightarrow +\infty$. 
We start by describing how to adapt the simultaneous breadth-first walk to the study of this phase.
We first obtain the concentration of the largest component and then study the fluctuations proving Theorem \ref{thm:barely-critical}.

\subsection{The concentration of the early giant component size}

In this setting we consider the initial configurations with $n$ vertices of mass equal to $1/n\epsilon_n$.
This scaling of mass is motivated by the fact that the early giant component is of size $O(n\epsilon_n)$ (see Section \ref{S:barely supercrit}). 
The random graph evolution is given by $(\mathcal{G}^{(n)}(t / (n \epsilon_n)^2 ))_{t \ge 0}$, so that the barely super-critical phase corresponds to the time window
\begin{equation}\label{eq:def_qnt}
	q_n(t) = (n \epsilon_n)^2 \frac{1 + t \epsilon_n}{n} = n \epsilon_n^2 (1 + t \epsilon_n), \ t>0.
\end{equation}
Therefore, the family of simultaneous breadth-first walks which correspond to our setting is
\begin{align*}
	Z_n^t (s) &:= \frac{1}{n \epsilon_n} \sum_{i = 1}^n \mathds{1}_{\{\bar{\xi}_i \cdot n \epsilon_n \le s q_n(t) \}} - s = \frac{1}{n \epsilon_n} \sum_{i = 1}^n \mathds{1}_{\{ \bar{\xi}_i  \le \epsilon_n (1 + t \epsilon_n) s \}} - s, \text{ for } t>0,
\end{align*}
where, as before, $(\bar{\xi}_i)_{i}$ are i.i.d.~exponentials with rate $1$.
Here we again rely on the results in \cite{Limic2019} already mentioned in Section \ref{S:sBFW_and_MC}.
In particular, the (ordered) list of excursion lengths of $(Z_n^t)_{t>0}$, viewed as a process in $t$, is equal in law to the process recording the (rescaled by $n \epsilon_n$ and ordered) connected component sizes of $\mathcal{G}^{(n)}\big( (1 + t \epsilon_n)/n \big)$ as $t>0$ increases.

Let us focus on the longest excursion of $Z_n^t$, which starts at some random time $S_n(t)/n \epsilon_n$ and ends at $(S_n(t) + L_n^{\mathrm{bsc}}(t))/n \epsilon_n$.
We again allow for a small abuse of notation, by identifying
the size of the largest component of $\mathcal{G}^{(n)}\big( (1 + t \epsilon_n)/n \big)$ 
with the length, multiplied by $n \epsilon_n$, of the longest excursion of $Z_n^t$.

Our approach is similar to that of Section \ref{S:proof_thm_enriquez}. The first difference from the super-critical case, is that the expectation of  $Z_n^t (s)$ depends on $n$.
Let us define
\[
\Upsilon_n^t(s) := \frac{1}{\epsilon_n}\mathbb{E}[Z_n^t (s)], \;\; \rho_n(t) := \rho\left(1+ \epsilon_n t\right),
\]
and note that 
\begin{equation}
	\label{expres:Upsilonnt}
	\Upsilon^t_n (s)  = \frac{1}{\epsilon_n^2} \left( 1 - \mathrm{e}^{- \epsilon_n (1 + t \epsilon_n) s}  - \epsilon_n s \right) = \frac{1}{\epsilon_n^2} \Phi^{1 + t \epsilon_n} (\epsilon_n s), \text{ for every } s \ge 0.
\end{equation}
Hence, $\Upsilon^t_n (\rho_n(t)/\epsilon_n) = 0$ for every $t>0$ and all $n\in \mathbb{N}$, and $\rho_n(t)/\epsilon_n$ is the unique zero of $\Upsilon^t_n$ in $(0,\infty)$. 
Furthermore, for each $t\geq 0$, $\Upsilon^t_n$ (resp.~$(\Upsilon^t_n)'(s) = \frac{1}{\epsilon_n} (\Phi^{1+t\epsilon_n})'(\epsilon_n s)$) converges to $\Upsilon^t: s \mapsto ts - s^2/2$ (resp.~$(\Upsilon^t)':s\mapsto t-s$) uniformly on compacts, as $n\to \infty$.
In particular, for each $t> 0$ the convergence of $\Upsilon^t_n$ to $\Upsilon^t$ as $n\to \infty$ yields 
\begin{equation}
	\label{eq:conv_rho_nt}
	\lim_{n\to \infty} \frac{\rho_n(t)}{\epsilon_n} = 2t,
\end{equation}
where the latter is the unique zero of $\Upsilon^t$ in $(0,\infty)$.

From now on we assume that $t>0$, and define 
\begin{equation}
	\label{D:X_bsc}
	X_n^t(s) :=  \sqrt{n \epsilon_n^3} \left( \frac{Z_n^t(s)}{\epsilon_n} - \Upsilon^t_n(s) \right), \text{ for every } s\geq 0.   
\end{equation}

The weak convergence of $(X_n^t)_n$ at a fixed $t>0$ can be obtained via a classical central limit theorem for martingales. It implies the following lemma: 
\begin{lemma}[Concentration and fluctuations of sBFW]\label{lemma:fluctuations_barely}
	For each $t>0$, as $n\to \infty$ 
	\[
	\left( \frac{Z_n^t(s)}{\epsilon_n} \right)_{s\geq 0} \Longrightarrow \Upsilon^t, \ \text{ and } \
	X_n^t\Longrightarrow  B,
	\]
	where $B$ is a standard Brownian motion, and where
	the weak convergence is in $\mathcal{C}([0,\infty))$. 
\end{lemma}
The proof of this lemma is postponed to Section \ref{sec:prooflem}.

Let us fix $t_1>0$.  
We again apply the Skorokhod's representation theorem. 
We claim there exist a family of random functions $ (\hat{X}_n^{t_1})_{n\geq 1}$ and a standard Brownian motion $(\hat{B}(t))_{t\geq 0}$  defined both on a same probability space, such that for each $n\geq 1$, $\hat{X}_n^{t_1}$ has the same distribution as  ${X}_n^{t_1}$, and for all $s_1>0$, 
\begin{equation}
	\label{equ:Z_fluctuates_bsc}
	\sup_{s \in [0, s_1]} \left| \hat{X}_n^{t_1}(s) - \hat{B}(s) \right|
	\xrightarrow[n \to \infty]{\text{a.s.}} 0.
\end{equation}
Set 
$\hat{Z}_n^{t_1}:= \epsilon_n( { \hat{X}_n^{t_1} }/{ \sqrt{n\epsilon^3_n} } +\Upsilon^{t_1}_n)$.
As in Section \ref{S:pf_thm_1}, from the equality in distribution  between $\hat{Z}_n^{t_1}$ and  $Z_n^{t_1}$ for each $n\geq 1$, we deduce that 
$$\hat{Z}_n^{t_1}= \frac{1}{n \epsilon_n} \sum_{i = 1}^n \mathds{1}_{\{ \bar{\xi}^{(n)}_i  \le \epsilon_n (1 + t \epsilon_n) s \}} - s,$$
where $(\bar{\xi}_i^{(n)})_{1\leq i\leq n}$ are a family of i.i.d. exponential rate 1 random variables. 
We now define $\hat{Z}_n^t$ for every $t>0$ as follows
\begin{align}
	\hat{Z}_n^t (s) & := \frac{1}{n \epsilon_n} \sum_{i = 1}^n \mathds{1}_{\{ \bar{\xi}^{(n)}_i  \le \epsilon_n (1 + t \epsilon_n) s \}} - s,   \text{ for }  s\geq 0.
\end{align}

The ``robustness'' of sBFW with respect to time-change is again very useful for extending the above asymptotics to the dynamical setting. Indeed, if we note that for each $t > 0$ 
\begin{equation}
	\label{equ:robust_bsc}
	\hat{Z}_n^{t}(s) = \hat{Z}_n^{t_1} \left( \frac{1 + t \epsilon_n}{1 + t_1 \epsilon_n} \cdot s \right) + \frac{t - t_1}{1 + t_1 \epsilon_n} \cdot \epsilon_n s,
\end{equation}
then the uniform convergence
$$
\sup_{s\in [0,s_1]} \left| \frac{\hat{Z}_n^{t_1}(s)}{\epsilon_n}  -\Upsilon^{t_1}(s) \right| \xrightarrow[n \to \infty]{\text{a.s.}} 0,  \text{ for every } s_1\in (0,\infty)
$$
immediately extends to the simultaneous (or doubly uniform) convergence
\begin{equation}
	\label{equ:doubly_uniform_convergence}
	\sup_{t\in [t_0,t_1]} \sup_{s\in [0,s_1]} \left| \frac{\hat{Z}_n^t(s)}{\epsilon_n}  - \Upsilon^t(s)  \right| \xrightarrow[n \to \infty]{\text{a.s.}} 0,  \ \forall s_1\in (0,\infty) \mbox{ and } \forall t_0 \in(0, t_1).
\end{equation}
Recall the ``past infimum'' operator $\operatorname{I}$, defined by $\operatorname{I}(f): s \mapsto \inf_{u\leq s} f(u)$.
Due to the Lipschitz continuity of $\operatorname{I}$ in the uniform topology, we immediately get that the almost-surely non-increasing process $\operatorname{I}(\hat{Z}_n^t/\epsilon_n)$ approximates 
$\operatorname{I}(\Upsilon^t)= \Upsilon^t(s) \mathds{1}_{\{s\geq 2t\}}$:
\begin{equation}
	\label{equ:doubly_uniform_convergence_I}
	\sup_{t\in [t_0,t_1]} \sup_{s\in [0,s_1]} \left| \operatorname{I}\left(\frac{\hat{Z}_n^t}{\epsilon_n}\right)(s)  - \operatorname{I}(\Upsilon^t)(s) \right| \xrightarrow[n \to \infty]{\text{a.s.}} 0,  \ \forall s_1\in (0,\infty) \mbox{ and } \forall t_0\in (0,t_1).
\end{equation}
Still, \eqref{equ:doubly_uniform_convergence}--\eqref{equ:doubly_uniform_convergence_I} together with additional properties of sBFW enable us to derive

\begin{lemma}[Control of the longest excursions in the barely super-critical regime]
	\label{lemma:control_in_bsc2}
	For every $t>0$, let $\hat{L}_n^{\mathrm{bsc}}(t)$ and $\hat{L}_n^{\mathrm{bsc};2}(t)$ denote the length of the longest and second-longest excursions  of  $\hat{Z}_n^{t}$, respectively. 
	For every $t_0 \in (0, t_1)$, 
	\begin{equation}\label{eq:conv_excursion_points}
		\sup_{t \in [t_0, t_1]} \left| \frac{\hat{L}_n^{\mathrm{bsc}}(t)}{n \epsilon_n} -  2t \right| \xrightarrow[n \to \infty]{\mathbb{P}} 0 \;\; \text{ and } \;\; \sup_{t \in [t_0, t_1]} \frac{\hat{L}_n^{\mathrm{bsc};2}(t)}{n \epsilon_n}  \xrightarrow[n \to \infty]{\mathbb{P}} 0.
	\end{equation}
\end{lemma}
\begin{remark}
	Notice that the point-wise version (for a fixed value of $t$) of \eqref{eq:conv_excursion_points} is well-known and already recalled in \eqref{eq:result_Nachmias_Peres}.
	We present a novel proof based on simultaneous breadth-first walks in Lemma \ref{lemma:control_in_bsc} in Section \ref{sec:prooflem}. 
	Moreover, using the fact that ${\hat{L}_n^{\mathrm{bsc}}(t)}/{n \epsilon_n}$ is increasing in $t$ and Dini's theorem one can easily get the first half of \eqref{eq:conv_excursion_points}.
	However, the uniform control on $t \in [t_0, t_1]$ of the length of the second longest excursion is more delicate, since  $(\hat{L}_n^{\mathrm{bsc}; 2}(t)/n \epsilon_n, \, t \geq t_0)$ is not monotone.
	The proof of this is also given in Section \ref{sec:prooflem}. 
\end{remark}
Thanks to Lemma \ref{lemma:control_in_bsc2} we can identify, analogously to the super-critical case, the longest excursion of $\hat{Z}_n^{t}$ as the only one of order $1$ occurring in $[0, 2t + \varepsilon't]$, for some $\varepsilon' > 0$.
As an immediate consequence of this fact and \eqref{equ:Z_fluctuates_bsc} we get that the process recording the starting point of the longest excursion, denoted by $(\hat{S}_n(t), t \ge t_0)$, satisfies
\begin{equation}
	\label{equ:unif_conv_S_bsc}
	\sup_{t \in [t_0,t_1]} \frac{\hat{S}_n(t)}{n \epsilon_n} \xrightarrow[n \to \infty]{\mathbb{P}}  0.
\end{equation}
\subsection{The fluctuations of the giant component size}
We now turn to the study of the second order fluctuations of $\hat{L}_n^{\text{bsc}}$.
Let us note
another immediate consequence of \eqref{equ:robust_bsc}: 
\[
\hat{X}_n^t(s) = \hat{X}_n^{t_1} \left( r_n(t) \cdot s \right), \text{ where } r_n(t)  = \frac{1 + t \epsilon_n}{1 + t_1 \epsilon_n}.
\]
This together with \eqref{equ:Z_fluctuates_bsc} yields that, for every $t\in [t_0,t_1]$:
\begin{align*}
	\sup_{s \in [0, s_1]} |\hat{X}_n^t(s) - \hat{B}(s)| &= \sup_{s \in [0, s_1]} |\hat{X}_n^{t_1}\big( r_n(t)s \big) - \hat{B}(s)| \\ 
	\le& \sup_{s \in [0, s_1]} \big| \hat{X}_n^{t_1}\big( r_n(t)s \big) - \hat{B} \big( r_n(t) s \big) \big| + \sup_{s \in [0, s_1]} \big| \hat{B}(s) - \hat{B} \big( r_n(t) s \big) \big|,
\end{align*}
implying that 
\begin{equation}
	\label{equ:doubly_unif_conv_X}
	\sup_{t \in [t_0, t_1]} \sup_{s \in [0, s_1]} |\hat{X}_n^t(s) - \hat{B}(s)| \xrightarrow[n \to \infty]{\text{a.s.}} 0,
\end{equation}   
due to \eqref{equ:Z_fluctuates_bsc}, the assumptions on $(\epsilon_n)_n$, and the uniform continuity of $\hat{B}$ on compact intervals.

Let us define
\begin{eqnarray}
	\label{D:sigma_bsc}
	\hat{\Sigma}_n^t &:=&\sqrt{n \epsilon_n^3} \Upsilon_n^t \left( \frac{\hat{S}_n(t) + \hat{L}_n^{\mathrm{bsc}}(t)}{n \epsilon_n} \right)\nonumber\\ 
	&\stackrel{\eqref{expres:Upsilonnt}}{=}&
	\sqrt{n \epsilon_n^3}
	\left(
	\Upsilon_n^t \left( \frac{\hat{S}_n(t) + \hat{L}_n^{\mathrm{bsc}}(t)}{n \epsilon_n}\right) -
	\Upsilon_n^t\left( \frac{\rho_n(t) }{\epsilon_n}\right)
	\right).
\end{eqnarray}
Applying the mean value theorem to $\Upsilon_n^t$ at $\rho_n(t)/\epsilon_n$ we get that
\begin{align}
	\label{expr_Sigma1}
	\hat{\Sigma}_n^t  = \big(\Upsilon_n^t\big)' (\gamma_n^t) \cdot \sqrt{n \epsilon_n^3} \left( \frac{\hat{S}_n(t) + \hat{L}_n^{\mathrm{bsc}}(t)}{n \epsilon_n} - \frac{\rho_n(t)}{\epsilon_n} \right) ,
\end{align}
for some random $\gamma_n^t$ in the interval  
\[[\rho_n(t)/{ \epsilon_n} \wedge (\hat{S}_n(t) + \hat{L}_n^{\mathrm{bsc}}(t))/{n \epsilon_n},   
\rho_n(t)/{ \epsilon_n} \vee (\hat{S}_n(t) + \hat{L}_n^{\mathrm{bsc}}(t))/{n \epsilon_n}].
\]
Recall \eqref{eq:conv_excursion_points} and \eqref{equ:unif_conv_S_bsc}.
It is easy to see that the second derivative of $\Upsilon_n^t$ is uniformly bounded on compacts, more precisely, that $ \sup_{t\in [t_0,t_1]} \|(\Upsilon_n^t)''\|_\infty = (1 + t_1 \epsilon_n)^2 < \infty$, uniformly in $n$,
for all ${t_0\in (0,t_1)}$.
Hence, for any $t_0\in(0,t_1)$  we can apply the mean value theorem, and recall the uniform convergence $(\Upsilon_n^t)'\to (\Upsilon^t)'$ to conclude that
\[
\sup_{t \in [t_0, t_1]} | \big(\Upsilon_n^t\big)' (\gamma_n^t) - \big(\Upsilon^t\big)' (2t) | \leq C
\sup_{t \in [t_0, t_1]} |  \gamma_n^t - 2t | +
\sup_{t \in [t_0, t_1]} | \big(\Upsilon_n^t\big)' (2t) - \big(\Upsilon^t \big)'(2t) |.
\]
Hence,
\begin{align}
	\label{equ:M_derivative_est}
	\sup_{t \in [t_0, t_1]} | \big(\Upsilon_n^t\big)' (\gamma_n^t) - \big(\Upsilon^t\big)' (2t) | 
	\xrightarrow[n \to \infty]{\text{a.s.}} 0.
\end{align}

Observe that $(\Upsilon^t)'(2t) = -t$.
Furthermore, recalling the definitions 
\eqref{D:X_bsc} and \eqref{D:sigma_bsc}
of $X_n^t$ and
$\Sigma_n^t$, respectively, 
we obtain an alternative expression for $\hat{\Sigma}_n^t$:
\begin{align}
	\label{expr_Sigma2}
	\hat{\Sigma}_n^t &= \sqrt{n \epsilon_n} \hat{Z}_n^t \left( \frac{\hat{S}_n(t) + \hat{L}_n^{\mathrm{bsc}}(t)}{n \epsilon_n} \right) -\hat{X}^{t}_{n}\left(\frac{\rho_n(t)}{\epsilon_n}\right)- \hat{R}_n^t,
\end{align}
where the error term is
\[
\hat{R}_n^t = \hat{X}_n^t \left( \frac{\hat{S}_n(t) + \hat{L}_n^{\mathrm{bsc}}(t)}{n \epsilon_n} \right) - \hat{X}_n^t \left( \frac{\rho_n(t)}{\epsilon_n} \right).
\]
As a consequence of \eqref{eq:conv_excursion_points}--\eqref{equ:doubly_unif_conv_X}, we have the following lemma:
\begin{lemma}[Error terms control]\label{lemma:error_terms_barely}
	For every $t_0\in(0,t_1)$, the non-negative processes
	\[ 
	\left( \sqrt{n \epsilon_n^3} \cdot \frac{\hat{S}_n(t)}{n \epsilon_n} \right)_{t \in [t_0, t_1]}, \;\;
	\left( - \sqrt{n \epsilon_n} \hat{Z}_n^t \left( \frac{\hat{S}_n(t) + \hat{L}_n^{\mathrm{bsc}}(t)}{n \epsilon_n} \right) \right)_{t \in [t_0, t_1]} \text{ and } \left( \hat{R}_n^t \right)_{t \in [t_0, t_1]}
	\]
	converge in probability uniformly to the zero process.
\end{lemma}
Before proving this lemma, let us end the proof of Theorem \ref{thm:barely-critical}.

Applying Lemma \ref{lemma:error_terms_barely} to the terms of formula \eqref{expr_Sigma2}, we obtain that the weak limit of $(\hat{\Sigma}_n^t)_{t \in [t_0, t_1]}$ as $n\to \infty$ is the weak limit of 
$\left( -\hat{X}_n^t \left( \frac{\rho_n(t)}{\epsilon_n} \right) \right)_{t \in [t_0, t_1]}$.
Due to Lemma \ref{lemma:fluctuations_barely}, the weak limit of $(\hat{\Sigma}_n^t)_{t \in [t_0, t_1]}$ is $(-\hat{B}(2t))_{t\in[t_0,t_1]}$. 

Let us isolate the term we are interested in by rewriting  equation \eqref{expr_Sigma1} slightly differently: 
\[
\big(\Upsilon_n^t\big)' (\gamma_n^t) \sqrt{n \epsilon_n^3} \left( \frac{\hat{L}^{\mathrm{bsc}}_n(t)}{n \epsilon_n} - \frac{\rho_n(t)}{\epsilon_n} \right) = \hat{\Sigma}_n^t - \big(\Upsilon_n^t\big)' (\gamma_n^t) \sqrt{n \epsilon_n^3}  \frac{\hat{S}_n(t) }{n \epsilon_n}.
\]
By \eqref{equ:M_derivative_est}, $(\big(\Upsilon_n^t\big)' (\gamma_n^t))_n$ converges to $-t$ in probability uniformly in $t\in[t_0,t_1]$, while  $\sqrt{n \epsilon_n^3}  \frac{\hat{S}_n(t) }{n \epsilon_n}$ converges to 0 in probability uniformly in $t\in[t_0,t_1]$  due to Lemma \ref{lemma:error_terms_barely}.

These observations combined imply that the process $\left(\sqrt{n \epsilon_n^3} \left( \frac{\hat{L}^{\mathrm{bsc}}_n(t)}{n \epsilon_n} - \frac{\rho_n(t)}{\epsilon_n} \right)\right)_{t\in[t_0,t_1]}$ converges weakly to $\left(\frac{1}{t}\hat{B}(2t)\right)_{t\in[t_0,t_1]}$. This ends the proof of Theorem \ref{thm:barely-critical} since  $(\hat{L}^{\mathrm{bsc}}_n(t))_{t \in[t_0,t_1]}$ has the same distribution as $(L^{\mathrm{bsc}}(t))_{t\in[t_0,t_1]}$ by Proposition 5 in \cite{Limic2019}.

\subsection{Proof of Lemmas \ref{lemma:fluctuations_barely}, \ref{lemma:control_in_bsc2} and \ref{lemma:error_terms_barely}\label{sec:prooflem}}
Let us first state some properties we will use in the proof of Lemma 4.1. Their proof is standard (see for example Theorem III.5.21 page 120 and Corollary III.5.2 in \cite{Protter2004} for the Doob-Meyer decomposition). For the convenience of readers, a proof is presented in Appendix \ref{sec:app_proof_lemma}.  
\begin{lemma}
	\label{lem:DoobMeyer_prelim}
	Let $\xi$ be a rate $q$ exponential random variable.
	Then, the process
	\[
	\Big( \mathds{1}_{\{\xi \le s\}} - q \cdot s \Big)_{s \ge 0}
	\]
	is a super-martingale with Doob--Meyer decomposition $M - A$ where $A$ is the increasing process defined by
	\(
	A(s) = { q} ( s - \xi)^+
	\)
	and $M$ is a martingale with predictable quadratic variation
	\(
	\langle M \rangle_s = q \min\{ s, \xi \}.
	\)
	Furthermore,
	\begin{itemize}
		\item[(i)] $\forall s\geq 0$, $\mathbb{E}[A(s)] = \mathrm{e}^{- q s} - 1 + q s,$ 
		\item[(ii)] $\forall \,  0\leq s_1<s_2$,
		\begin{equation*}
			\operatorname{Var}[A(s_2)-A(s_1)]=2 \frac{\left(1-\mathrm{e}^{-q (s_2 - s_1) }-q (s_2 - s_1) \mathrm{e}^{-q (s_2 - s_1) }\right)}{\mathrm{e}^{qs_1}} - \left( \frac{1-\mathrm{e}^{-q (s_2 - s_1)} }{\mathrm{e}^{q s_1}}\right)^2.
		\end{equation*}
	\end{itemize}
\end{lemma}

\begin{proof}[Proof of Lemma \ref{lemma:fluctuations_barely}]
	Recall the definition \eqref{D:X_bsc} of $X^{t}_n$.
	Let us find the semi-martingale decomposition of $Z_n^t(s) / \epsilon_n $ by using the Lemma \ref{lem:DoobMeyer_prelim}.
	We have
	\begin{align}
		\frac{Z_n^t(s)}{\epsilon_n} &= \frac{1}{n \epsilon_n^2} \sum_{i = 1}^n \mathds{1}_{\left\{ \overline{\xi}_i \le \epsilon_n (1 + t \epsilon_n) s \right\}} - \frac{s}{\epsilon_n} \nonumber \\
		&= \frac{1}{n \epsilon_n^2} \sum_{i = 1}^n \left( \mathds{1}_{\left\{ \overline{\xi}_i \le \epsilon_n (1 + t \epsilon_n) s \right\}} - \epsilon_n (1 + t \epsilon_n)s \right) + \frac{\epsilon_n (1 + t \epsilon_n)s}{\epsilon_n^2}  -  \frac{s}{\epsilon_n} \nonumber \\
		&= \frac{1}{n \epsilon_n^2} \sum_{i = 1}^n M_{i,n}(s) - \frac{1}{n \epsilon_n^2} \sum_{i = 1}^n A_{i,n}(s) + ts, \label{eq:DMdecomposition}
	\end{align}
	where
	\[
	A_{i,n}(s) := \epsilon_n (1 + t \epsilon_n) \left(\, s - \frac{\overline{\xi}_i}{\epsilon_n (1 + t \epsilon_n)} \right)^+
	\]
	and
	\[
	M_{i,n}(s) := A_{i,n}(s)+\mathds{1}_{\left\{ \overline{\xi}_i \le \epsilon_n (1 + t \epsilon_n) s \right\}} - \epsilon_n (1 + t \epsilon_n) s.
	\]
	The first term in the previous sum is the martingale part and the other two terms form the finite variation part.
	Let us denote
	\[
	M_n^t (s) :=  \frac{1}{n \epsilon_n^2} \sum_{i = 1}^n M_{i,n}(s) \text{ and }  A_n^t (s) :=  \frac{1}{n \epsilon_n^2} \sum_{i = 1}^n A_{i,n}(s).
	\]
	It is clear that $\mathbb{E}[A_n^t] = st-\mathbb{E}[Z_n^t(s) / \epsilon_n]$.
	Hence,
	\begin{equation}
		\label{eq:X_n_decomposition}
		X_n^t(s) = \sqrt{n \epsilon_n^3} M_n^t(s) + \sqrt{n \epsilon_n^3} \left( A_n^t(s) - \mathbb{E}[A_n^t (s)] \right).
	\end{equation}
	We will study separately the convergence of the martingale and the finite variation parts.
	Let us first consider the finite variation part $$\Delta_n(s):=\sqrt{n \epsilon_n^3} \left( A_n^t(s) - \mathbb{E}[A_n^t (s)] \right) =\frac{1}{\sqrt{n\epsilon_n}}\sum_{i=1}^{n}\left( A_{i,n}(s) - \mathbb{E}[A_{i,n} (s)] \right).$$
	Take two reals $0\leq s_1<s_2$ and set $\alpha_n^t=\epsilon_n(1+t\epsilon_n)$.
	
	By Lemma \ref{lem:DoobMeyer_prelim}, 
	\begin{eqnarray*}
		\mathbb{E}[ (\Delta_n(s_2)-\Delta_n(s_1))^2 ] &=&\frac{1}{\epsilon_n}\operatorname{Var}\Big(  A_{1,n}(s_2)-A_{1,n}(s_1) \Big) \\
		&=&\frac{2\mathrm{e}^{-\alpha^{t}_n s_1}}{\epsilon_n}\left( 1-\mathrm{e}^{-\alpha^{t}_n (s_2-s_1)}-\alpha^{t}_n (s_2-s_1) \mathrm{e}^{-\alpha^{t}_n (s_2-s_1)}\right)\\
		&&-\frac{ \mathrm{e}^{-2\alpha^{t}_n s_1}}{\epsilon_n}\left(1-\mathrm{e}^{-\alpha_n^t(s_2-s_1)}\right)^2.
	\end{eqnarray*}
	Using that for all $x \geq 0$, we have $1 - \mathrm{e}^{-x}-x \mathrm{e}^{-x}\leq x^2$, we obtain
	\begin{equation}
		\mathbb{E}((\Delta_n(s_2)-\Delta_n(s_1))^2)\leq 2\epsilon_n(1+t\epsilon_n)^2(s_2-s_1)^2.   
	\end{equation}
	This upper bound implies that for every $s \ge 0$, the sequence $(\Delta_n(s))_n$ converges to $0$ in $L^2$, as well as the tightness of $(\Delta_n)_n$ by Billingsley's criterion (\cite{Billingsley}, Theorem 13.5) adapted for the case of processes indexed by $\mathbb{R}_+$ (\cite{Jacod}, Theorem 4.1). 
	Therefore, $(\Delta_n)_n$ converges weakly to the zero process.

	For the martingale part, note that its quadratic variation satisfies
	\begin{align*}
		\left[ \sqrt{n \epsilon_n^3} M_n^t \right]_s &= \left[ \frac{\sqrt{n \epsilon_n^3}}{\epsilon_n} Z_n^t \right]_s = \frac{1}{n \epsilon_n} \sum_{i = 1}^n \mathds{1}_{\{ \bar{\xi_i} \le \epsilon_n (1 + t \epsilon_n)s \} } \xrightarrow[n \to \infty]{ \mathbb{P} } s,
	\end{align*}
	by Chebyshev inequality, since the variance of the sum in the above identity is $O(1/(\epsilon_n n))$.
	
	Since any jump of $( \sqrt{n \epsilon_n^3} M_n^t )$ is of size $\sqrt{n \epsilon_n^3}/(n \epsilon_n^2) = 1/\sqrt{n \epsilon_n} \to_{n\to \infty} 0$, the martingale functional central limit theorem \cite[Theorem 7.1.4, p.\ 339]{EthierKurtz} implies that $( \sqrt{n \epsilon_n^3} M_n^t )_n$ converges weakly to the standard Brownian motion.
	
	Recall \eqref{eq:X_n_decomposition}.
	Since we showed above that $(\Delta_n)_n$ converges to the zero process,
	we can conclude that $(X_n^t)_n$ also converges weakly to the standard Brownian motion.
	In particular,  $(Z_n^t/\epsilon_n)_n$ converges to $\lim_n \Upsilon_n^t=\Upsilon^t$ as $n\to \infty$.
\end{proof}
In order to prove Lemma \ref{lemma:control_in_bsc2} we need the following auxiliary result.
\begin{lemma}[Static control of the longest excursions in the barely super-critical regime]
	\label{lemma:control_in_bsc}
	For every $t \in [t_0,t_1]$, 
	\begin{equation*}\label{eq:1st_2nd_CC_bsc}
		\left| \frac{\hat{L}_n^{\text{bsc}}(t)}{n \epsilon_n} -  2t \right| \xrightarrow[n \to \infty]{\mathbb{P}} 0 \;\; \text{ and } \;\;  \frac{\hat{L}_n^{\text{bsc};2}(t)}{n \epsilon_n}  \xrightarrow[n \to \infty]{\mathbb{P}} 0,
	\end{equation*}
\end{lemma}

\begin{remark}
	This result in terms of the sizes of the longest and the second longest component of $\mathcal{G}^{(n)}((1+t \epsilon)/n)$ is well-known, as we commented in \eqref{eq:result_Nachmias_Peres}.
	Still, for the sake of completeness, we provide another proof based on the sBFW.
\end{remark}

\begin{proof}[Proof of Lemma \ref{lemma:control_in_bsc}]
	Note that the unique excursion (above past minimum) of the limit $(\Upsilon^t(s),\,s \geq 0)$ of $\hat{Z}_n^t/\epsilon_n$ is $[0,2t]$. By the argument relying on \eqref{equ:tall_excursion_f_implies_excursion_g}--\eqref{equ:decrease_f_implies_no_excursion_g} from Section \ref{S:concentration_sc} and property \eqref{equ:doubly_uniform_convergence} we can conclude that
	as $n\to \infty$, on the events
	\[
	\left\{ \sup_{s\in [0, 6t]} \left| \frac{1}{\epsilon_n} \hat{Z}_n^t(s) - \Upsilon^t(s) \right| < \frac{1}{(n \epsilon^{3}_n)^{1/4}} \right\} \text{ which occur w.h.p.}
	\]
	we have both
	\begin{enumerate}
		\item[(i)] $\hat{Z}_n^t$ has an excursion of length approximately $2t$ which begins at 
		\[
		{\bar{\xi}_{(I_n^t)}^{(n)}/(\epsilon_n(1+t\epsilon_n))= o_n(1)},
		\]
		for some random index $I_n^t$; and
		\item[(ii)] $\hat{Z}_n^t$ has no other excursion with length of order $1$  which starts in 
		$[0, 5t]$.
	\end{enumerate}
	Due to these properties, in order to prove Lemma \ref{lemma:control_in_bsc} it remains to show that the longest excursion of $\hat{Z}^t_n$ in $[0,5 t]$ is indeed its longest excursion in $\mathbb{R}_+$, and that the second longest excursion is $o_{\mathbb{P}}(1)$.
	
	It is well-known that excursions of $\hat{Z}_n^t$ appear in size-biased order with respect to their lengths.
	The following related fact, particularly useful for our conclusion, was recently proved in \cite[Lemma 3.12]{2024blanc}: given an ordered list $L_1, L_2, \ldots$ of excursion lengths of $\hat{Z}_n^t$, which respectively begin at times $T_1, T_2, \ldots$, the random vector of ``local times''  ordered according to the excursion length list given above
	\begin{equation*}
		\Big( -\hat{Z}_n^t(T_i-)= -\inf_{s\leq T_i} Z_n^t(s) = -\operatorname{I}(\hat{Z}_n^t)(T_i)
		\Big)_{i \ge 1}
	\end{equation*}
	is conditionally distributed, given $L_1,L_2,\ldots$, as a family of independent exponential random variables, where 
	$-\hat{Z}_n^t(T_k-)$ has exponential (rate $q_n(t) \cdot L_k$) distribution, and where $q_n(t)$ is given by \eqref{eq:def_qnt}.
	In particular, on the event that two longest excursions of $\hat{Z}_n^t$ have lengths $L_1^n(t)$ and $L_2^n(t)$, the conditional distribution of \[-\operatorname{I}(\hat{Z}_n^t/\epsilon_n)(T_1^n(t))=-\hat{Z}_n^t(T_1^n(t)-)/\epsilon_n\text{ and } -\operatorname{I}(\hat{Z}_n^t/\epsilon_n)(T_2^n(t))=-\hat{Z}_n^t(T_2^n(t)-)/\epsilon_n\] is that of two independent exponentials with rates 
	\[
	n\epsilon_n^3(1+t\epsilon_n) L_1^n(t) \text{ and } 
	{n\epsilon_n^3(1+t\epsilon_n) L_2^n(t)},
	\]
	respectively.
	
	For the sake of simpler notation let us abbreviate $\operatorname{I}_n^t \equiv \operatorname{I}(\hat{Z}_n^t/\epsilon_n)$.
	The above reasoning implies that, for any fixed $\delta\in(0,1/2)$ we get 
	\begin{align*}
		\mathbb{P} &\left( L_1^n(t)\geq  L_2^n(t) >\delta t , (-\operatorname{I}_n^t(T_1^n(t))) \vee (-\operatorname{I}_n^t(T_2^n(t))) \ge \frac{1}{\sqrt{n \epsilon_n^3}} \right) \\
		&\le \mathbb{P}\left( L_1^n(t)\geq  L_2^n(t) >\delta t , \left( -\operatorname{I}_n^t(T_1^n(t)) \cdot L_1^n(t)\right) \vee \left(-\operatorname{I}_n^t(T_2^n(t)) \cdot L_2^n(t) \right) \ge \frac{ \delta t }{\sqrt{n \epsilon_n^3}} \right)\\
		&\le \mathbb{P}\left( \vartheta_1 \vee \vartheta_2 > \delta t \sqrt{n \epsilon_n^3} (1 + t \epsilon_n) \right) \xrightarrow[n \to \infty]{} 0,
	\end{align*}
	where $\vartheta_1$ and $\vartheta_2$ are independent exponential rate $1$ random variables, and where we used the hypotheses on $(\epsilon_n)_n$.
	
	In words, on the event that the longest two excursions of $\hat{Z}_n^t$ are both longer than $\delta t$, both values $\operatorname{I}_n^t(T_1^n(t))=\hat{Z}_n^t(T_1^n(t)-)/\epsilon_n$  
	and
	$\operatorname{I}_n^t(T_2^n(t))=\hat{Z}_n^t(T_2^n(t)-)/\epsilon_n$
	converge to $0$ in probability as $n\to \infty$.
	
	Recalling \eqref{equ:doubly_uniform_convergence_I} and noting $\operatorname{I}(\Upsilon^t)((2+\delta)t) =  - \delta(1 + \delta/2) t^2$,
	implies that 
	\[
	\sup_{s\geq (2+\delta) t} \operatorname{I}_n^t(s) = \operatorname{I}_n^t \big( (2+\delta)t \big) \leq - \delta t^2
	\]
	with overwhelming probability as $n\to \infty$.
	We have argued that the event 
	\[
	{\{L_1^n(t) \geq  L_2^n(t) >\delta t\}}
	\] 
	has asymptotically the same probability as
	\[
	\left\{ L_1^n(t) \geq  L_2^n(t) >\delta t, \operatorname{I}_n^t((2+\delta)t)\leq - \delta t^2, (-\operatorname{I}_n^t(T_1^n(t))) \vee (-\operatorname{I}_n^t(T_2^n(t))) \le \frac{1}{\sqrt{n \epsilon_n^3}} \right\}
	\]
	and this event is included for all large $n$ in the event
	\begin{align}
		\label{E:above event}
		\{L_1^n(t) \geq  L_2^n(t) >\delta t, \, T_1^n(t) \vee T_2^n(t)\leq (2+ \delta) t\}.
	\end{align}
	
	On the above event there are (at least) two excursions of $\hat{Z}_n^t$ of length greater than $\delta t$ 
	starting before $2.5 t$ (recall that $\delta<1/2$). Recalling (ii) from the initial paragraph of this proof, the occurrence of this event has negligible probability as $n\to \infty$.
	We therefore derived that for any $\delta>0$, $\{ \hat{L}_n^{\text{bsc;2}}(t)\leq \delta t\}$ occurs with overwhelming probability as $n\to \infty$.
	Moreover, since we already know (see (i) above) that w.h.p.~there exists an excursion of $\hat{Z}_n^t$ of an approximate length $2t$ in $[0,5t]$ (in fact, already in $[0,(2+\epsilon)t]$), this excursion must also be its globally longest excursion w.h.p.
\end{proof}

Let us note that the statement of Lemma \ref{lemma:control_in_bsc} also holds on the initial probability space for the sequences $(L^{\mathrm{bsc}}_n)_n$ and $(L^{\mathrm{bsc};2}_n)_n$ since $(\hat{Z}_n^t)_{t\in [t_0,t_1]}$ has the same distribution as  $(Z_n^t)_{t \in [t_0,t_1]}$. Likewise, instead of proving Lemma \ref{lemma:control_in_bsc2} we can prove the following  equivalent version on the initial space. 
\begin{lemma}
	\label{lemma:control_in_bsc2b}
	For every $t_0 \in (0, t_1)$, 
	\begin{equation}\label{eq:conv_excursion_points2}
		\sup_{t \in [t_0, t_1]} \left| \frac{{L}_n^{\mathrm{bsc}}(t)}{n \epsilon_n} -  2t \right| \xrightarrow[n \to \infty]{\mathbb{P}} 0 \;\; \text{ and } \;\; \sup_{t \in [t_0, t_1]} \frac{{L}_n^{\mathrm{bsc};2}(t)}{n \epsilon_n}  \xrightarrow[n \to \infty]{\mathbb{P}} 0.
	\end{equation}
\end{lemma}
While the proof of Lemma \ref{lemma:control_in_bsc} relied on the 
sBFW, we shall prove Lemma \ref{lemma:control_in_bsc2b} by arguments using directly the random graph evolution.

\begin{proof}[Proof of Lemma \ref{lemma:control_in_bsc2}]
	Let us denote by $\iota_n : s \mapsto (s n - 1)/ \epsilon_n$, the inverse of the time change function $s_n: t \mapsto (1 + \epsilon_n t)/n$.
	Fix some $\delta\in (0, t_0)$ and for each $n$ define two  stopping times, $T_n\equiv T_n(\delta;t_0)$ and $V_n\equiv V_n(\delta;t_0)$, with respect to the natural filtration associated with the $n$th random graph process $\mathcal{G}^{(n)}$:
	\begin{eqnarray*}
		T_n&:=&\inf\left\{s \geq \frac{1 + t_0 \epsilon_n}{n}: \frac{L_n^{\mathrm{bsc}; 2}\big( \iota_n(s) \big)}{\epsilon_n n} > 5\delta\right\}, \\
		V_n&:=&\inf\left\{ s \geq \frac{1 + t_0 \epsilon_n}{n}: \left|\frac{L_n^{\mathrm{bsc}}(\iota_n(s))}{\epsilon_n n} - 2 \iota_n(s) \right| > \delta\right\}.
	\end{eqnarray*}
	
	By Lemma  \ref{lemma:control_in_bsc} and Dini's theorem we get the uniform convergence on compacts of $L_n^{\mathrm{bsc}}/(2\epsilon_n n)$ to the identity function, which corresponds to the first half of 
	\eqref{eq:conv_excursion_points2}.
	Hence, we have that for every $C \in (0, +\infty)$:
	\begin{equation}
		\label{eq:V_is_infinite_bsc}
		\lim_{n\to \infty} \mathbb{P} \left(V_n < \frac{1 + C \epsilon_n}{n} \right) =0.
	\end{equation}
	Thus, to conclude the proof it suffices to show that $\mathbb{P}(T_n \wedge V_n \leq (1 + t_1 \epsilon_n)/n)$ is asymptotically negligible.
	
	Due to the (strong) Markov property of $\mathcal{G}^{(n)}$
	the multiplicative coalescent evolution of
	the component masses in $$
	\left( \mathcal{G'}^{(n)}(u), u \geq 0 \right) :=
	\left(\mathcal{G}^{(n)}\left(V_n \wedge T_n \wedge \frac{1 + t_1 \epsilon_n}{n} + u \right), \, u\geq 0\right),$$
	given 
	$\mathcal{G}^{(n)}(V_n \wedge T_n \wedge \frac{1 + t_1 \epsilon_n}{n})$, is
	conditionally independent 
	of \[(\mathcal{G}^{(n)}(u),\, u < V_n\wedge T_n \wedge \frac{1+ t_1 \epsilon_n}{n}).\]
	
	Given $\mathcal{G'}^{(n)}(0)=\mathcal{G}^{(n)}(V_n \wedge T_n \wedge \frac{1 + t_1 \epsilon_n}{n})$,
	we will focus on the two components of $\mathcal{G'}^{(n)}$ which respectively contain the largest and the second largest  component of $\mathcal{G'}^{(n)}(0)$ which merge at a new stopping time $V_n'$. Note that $(T_n\wedge V_n\wedge \frac{1 + t_1 \epsilon_n}{n}) + V_n'$ is a stopping time with respect to  $\mathcal{G}^{(n)}$. 
	
	On the event 
	$\{T_n < V_n \wedge \frac{1 + t_1 \epsilon_n}{n} \}$, we know that 
	the largest (resp.~second largest) component of $\mathcal{G'}^{(n)}(0) $ has size greater than or equal to $$( 2 \iota_n( T_n ) -\delta) n \epsilon_n \geq (2 t_0 - \delta) n \epsilon_n \geq t_0 n \epsilon_n\quad (\text{resp.~} 5\delta n \epsilon_n ).$$
	Hence, the merger at time $V_n'$
	yields a component of size greater than or equal to $(2 {\iota_n( T_n )} + 4 \delta )n \epsilon_n$ in $\mathcal{G'}^{(n)}(V_n') = \mathcal{G}^{(n)} (V_n \wedge T_n \wedge \frac{1 + t_1 \epsilon_n}{n} + V_n')$.
	On $\{T_n < V_n \wedge \frac{1 + t_1 \epsilon_n}{n} \}$, the conditional probability of $\{V_n'> \delta \epsilon_n/n\}$
	is bounded above by
	\begin{equation}
		\label{eq:Vnprime_cannot_be_large}
		\exp\left\{ -\frac{\delta \epsilon_n}{n} \cdot \left(2 {\iota_n(T_n)}  -\delta \right) n \epsilon_n  \cdot 5 \delta n \epsilon_n \right\} \le \exp \left\{ -5 t_0 \delta^2   n \epsilon_n^3 \right\}.
	\end{equation}
	Due to our hypothesis on $(\epsilon_n)_n$ this upper bound converges to $0$ as $n\to \infty$. As already discussed, the event
	$$ 
	E_n^{\le} := \left\{T_n < V_n \wedge \frac{1 + t_1 \epsilon_n}{n},\,
	V_n' \leq \frac{\delta \epsilon_n}{n} \right\}
	$$
	is included in
	$$
	\left\{ 
	L_n^{\mathrm{bsc}}(\iota_n(T_n) + \delta ) \geq 
	\left( 2 {\iota_n(T_n)}  + 4 \delta \right) n \epsilon_n
	\right\}.
	$$ 
	Using the fact that $\iota_n(T_n)+\delta= \iota_n(T_n + \delta \epsilon_n/n )$, the above implies
	\begin{align*}
		E_n^{\le} &\subset \left\{\frac{L_n^{\mathrm{bsc}}(\iota_n(T_n) + \delta)}{ n \epsilon_n } \geq  2\left( {\iota_n(T_n)} + \delta \right)  + 2\delta, T_n < \frac{1 + t_1 \epsilon_n}{n} \right\} \\
		&\subset \left\{V_n \leq T_n + \frac{\delta \epsilon_n}{ n}, T_n < \frac{1 + t_1 \epsilon_n}{n} \right\} \\
		&\subset \left\{V_n \leq \frac{1 + (t_1 + \delta) \epsilon_n}{ n} \right\}.
	\end{align*}
	Define
	$$ 
	E_n^{>} := \left\{T_n < V_n \wedge \frac{1 + t_1 \epsilon_n}{n},\,
	V_n' > \frac{\delta \epsilon_n}{n} \right\}.
	$$
	Combining all of the above we can now conclude
	\begin{align}
		\nonumber \mathbb{P}\left( T_n \wedge V_n \leq \frac{1 + t_1 \epsilon_n}{n} \right) \leq & \, \mathbb{P}\left( V_n \leq T_n \wedge \frac{1 + t_1 \epsilon_n}{n}  \right) + \mathbb{P}(E_n^\le) + \mathbb{P}(E_n^>) \\
		\leq&
		\label{eq:bound_tail_bsc} \,
		2 \, \mathbb{P} \left( V_n \leq \frac{1 + (t_1 + \delta) \epsilon_n}{n} \right) + \mathbb{P}(E_n^>). 
	\end{align}
	The second term in \eqref{eq:bound_tail_bsc} decays to zero with $n$ diverging, as was already indicated 
	below \eqref{eq:Vnprime_cannot_be_large}. 
	Due to 
	\eqref{eq:V_is_infinite_bsc}, the first term in \eqref{eq:bound_tail_bsc} is also negligible as $n\to \infty$.
\end{proof}

\begin{proof}[Proof of Lemma \ref{lemma:error_terms_barely}]
	Recall that
	\[
	\hat{R}_n^t =  \hat{X}_n^t \left( \frac{ \hat{S}_n(t) +  \hat{L}_n^{\mathrm{bsc}}(t)}{n \epsilon_n} \right) -  \hat{X}_n^t \left( \frac{\rho_n(t)}{ \epsilon_n} \right).
	\]
	As a consequence of 
	\eqref{eq:conv_rho_nt},
	\eqref{eq:conv_excursion_points2}, \eqref{equ:unif_conv_S_bsc} and 
	\eqref{equ:doubly_unif_conv_X},
	it is clear that $( \hat{R}_n^t)_{t \in [t_0, t_1]}$ 
	converges  uniformly to the zero process. 
	
	The rest of the argument is analogous to the proof of Lemma \ref{lemma:control}.  
	In particular, we  note that both processes 
	$$\left(\sqrt{n \epsilon_n^3} \cdot \frac{ \hat{S}_n(t)}{n \epsilon_n} \right)_{t \in [t_0, t_1]} \;\; \text{ and } \;\; \left(- \sqrt{n \epsilon_n}  \hat{Z}_n^t \left( \frac{ \hat{S}_n(t) +  \hat{L}_n^{\mathrm{bsc}}(t)}{n \epsilon_n} \right) \right)_{t \in [t_0, t_1]}$$ 
	are decreasing in $t$ for all $n$ large enough, so it suffices to show the convergence for a fixed $t$.
	
	So let us fix some $t \in [t_0,t_1]$, and consider the problem $ \hat{S}_n(t) = o_{\mathbb{P}}(\sqrt{n/\epsilon_n})$ as $n\to \infty$.
	The idea is (compare with Lemma \ref{lemma:control})
	to check that $ \hat{Z}_n^t$ is positive on a segment $[\delta /\sqrt{n \epsilon_n^3},\alpha(t)]$, with probability converging to $1$ as $n\to \infty$, where $\alpha(t)>0$ does not depend on $n$. 
	The divergence of $(n\epsilon_n^3)_n$ is again important.
	
	The uniform convergence on compact sets of $(\Upsilon_n^{t})'$ to $(\Upsilon^{t})':s\mapsto t-s$ implies that:
	\begin{itemize}
		\item we can choose $\alpha(t)>0$ such that $\Upsilon_n^t$ is increasing on $[0, \alpha(t)]$ for all large $n$,
		\item  if $(a_n)_n$ is a positive sequence that tends to $+\infty$, then
		\begin{equation}
			\label{eqn:propdervUpsilon}
			\left(a_n \Upsilon_n^t\left(1/a_n\right) \right)_n \text{ converges to }   t.
		\end{equation}
	\end{itemize}
	By \eqref{equ:Z_fluctuates_bsc}, 
	\begin{equation}
		\label{equ:conseqence_diff_approx_Z}
		\sup_{s \in K} |\sqrt{n \epsilon_n} \hat{Z}_n^t(s) - \sqrt{n \epsilon_n^3} \Upsilon_n^t(s) - \hat{B}(s) |  = o_{\mathbb{P}}(1), \text{ on any compact set }K.
	\end{equation}
	Now consider $\delta > 0$ arbitrary small.
	
	We again split the interval $[ \delta/\sqrt{n \epsilon_n^3}, \alpha(t)]$ into two pieces.
	The continuity of $\hat{B}$, joint with  \eqref{eqn:propdervUpsilon} and  \eqref{equ:conseqence_diff_approx_Z}, yields  
	\[
	\sqrt{n \epsilon_n} \hat{Z}_n^t(s) \ge \sqrt{n \epsilon_n^3} \Upsilon_n^t \left( \frac{\delta}{\sqrt{n \epsilon_n^3}} \right) + o_{\mathbb{P}}(1), \ \ s \in [\delta/\sqrt{n \epsilon_n^3}, 1/ \sqrt[4]{n \epsilon_n^3}], 
	\]
	where the RHS above is larger than  $t \delta/2 >0$, for all $n$ sufficiently large.
	Similarly, the boundedness of $\hat{B}$ on $[1/ \sqrt[4]{n \epsilon_n^3}, \alpha(t)]\subset [0,\alpha(t)]$ and \eqref{equ:conseqence_diff_approx_Z} give
	\[
	\sqrt{n \epsilon_n} \hat{Z}_n^t(s) \ge \sqrt{n \epsilon_n^3 } \Upsilon_n^t \left( \frac{1}{\sqrt[4]{n \epsilon_n^3}} \right) + O_{\mathbb{P}}(1), \ \ 
	s \in [1/ \sqrt[4]{n \epsilon_n^3}, \alpha(t)].
	\]
	The property \eqref{eqn:propdervUpsilon} now implies the divergence of 
	$\sqrt{n  \epsilon_n^3 } \Upsilon_n^t \left( {1}/{\sqrt[4]{n \epsilon_n^3}}\right)$ to $+\infty$ and therefore  the positivity of $\hat{Z}_n^t$ on $[1/ \sqrt[4]{n \epsilon_n^3}, \alpha(t)]$ for all large $n$. 
	As argued above (and in Lemma \ref{lemma:control}) this concludes the proof of the first half of the lemma.
	
	As in the super-critical setting, by the properties of $\hat{Z}_n^t$, we have
	\[
	-\hat{Z}_n^t \left( \frac{\hat{S}_n(t) + \hat{L}_n^{\mathrm{bsc}}(t)}{n \epsilon_n} \right) = 
	-\hat{Z}_n^t \left( \frac{\hat{S}_n(t) }{n \epsilon_n} -\right) \leq \frac{\hat{S}_n(t) }{n \epsilon_n},
	\]
	but multiplying the inequality by $\sqrt{\epsilon_n n}$ does not give a sufficiently small scale on the right-hand side to be able to conclude the second half of the lemma directly from the first half, as done in the proof of Lemma \ref{lemma:control}.
	
	We prove that $\sqrt{n \epsilon_n} \hat{Z}_n^c(\hat{S}_n(c)/n \epsilon_n -)\to 0$ in probability as $n \to \infty$ as follows. 
	Fix some $\vartheta>0$, define $\mathcal{A}_n\equiv \mathcal{A}_n^t(\vartheta) := \{ \sqrt{n \epsilon_n^3} \frac{\hat{S}_n(t)}{n \epsilon_n} < \vartheta \}$, and 
	note that $\mathbb{P}(\mathcal{A}_n) \to 1$.
	On $\mathcal{A}_n$ we have
	\[
	0 \le - \hat{Z}_n^t \left( \frac{\hat{S}_n(t)}{n \epsilon_n} - \right) \le \sup_{s \in [0, \vartheta/\sqrt{n \epsilon_n^3} ] } - \hat{Z}_n^t(s) = -\operatorname{I}(\hat{Z}_n^t)\left(\frac{\vartheta}{\sqrt{n \epsilon_n^3}}\right).
	\]
	Furthermore, due to  \eqref{equ:conseqence_diff_approx_Z}, the continuity and superadditivity of the past infimum operator $\operatorname{I}$ we have 
	$  \operatorname{I}(\sqrt{n\epsilon_n} \hat{Z}_n^t) = 
	\operatorname{I}( \sqrt{n\epsilon_n^3}(\Upsilon_n^t) +\hat{B}) + o_{\mathbb{P}}(1)\geq 
	\sqrt{n\epsilon_n^3}\operatorname{I}(\Upsilon_n^t) + \operatorname{I}(\hat{B}) + o_{\mathbb{P}}(1)
	$ on any compact $K\subset\mathbb{R}_+$. 
	Evaluating the above inequality at $\vartheta/\sqrt{n\epsilon_n^3}$,
	and using  the fact that $\Upsilon_n^t$ is increasing on $[0,\alpha(t)]$ together with the continuity of $\hat{B}$
	we get  that on $\mathcal{A}_n$ 
	\[
	\sqrt{n \epsilon_n} \sup_{s \in [0, \vartheta/\sqrt{n \epsilon_n^3} ] } - \hat{Z}_n^t(s) \le - \sqrt{n\epsilon_n^3}\operatorname{I}(\Upsilon_n^t)\left( \frac{\vartheta}{\sqrt{n \epsilon_n^3}} \right) +  o_{\mathbb{P}}(1) = o_{\mathbb{P}}(1).
	\]
	As already argued, this yields the proof of the lemma.
\end{proof}

\appendix

\section{Proof of Lemma \ref{lem:DoobMeyer_prelim}}\label{sec:app_proof_lemma}

Let $\big( N(t),\, t \ge 0 \big)$ be a rate $q$ homogeneous Poisson process.
It is well-known that $\big( N(t) - q \cdot t,\, t \ge 0 \big)$, is a martingale with respect to its natural filtration, which is the same as the filtration generated by $\big( N(t),\, t \ge 0 \big)$.
Moreover, its quadratic variation and predictable quadratic variation at time $t \ge 0$ are $N(t)$ and $q \cdot t$, respectively.
Let us define the stopping time $\xi := \inf\{t \geq 0: N(t) \ge 1\}$ recording the first jump of $\big( N(t),\, t \ge 0 \big)$.
Notice that the martingale $\big( N(t) - q \cdot t,\, t \ge 0 \big)$ stopped at time $\xi$ is $(\mathds{1}_{\{\xi \le t\}} - q  (\xi \wedge t),\, t \ge 0)$, whose quadratic variation and predictable quadratic variation are $(\mathds{1}_{\{\xi \le t\}},\, t \ge 0)$ and  $(q (t \wedge \xi),\, t \ge 0)$, respectively.
Hence,
$
\big( \mathds{1}_{\{\xi \le t\}} - q \cdot t ,\, t \ge 0 \big)
$
is a super-martingale with Doob-Meyer decomposition $M-A$ where
$
{M(t) = \mathds{1}_{\{\xi \le t\}} - q (\xi \wedge t)}
$, $[ M ]_t = \mathds{1}_{\{\xi \le t \}}$, $\langle M \rangle_t =q(\xi \wedge t)$ and ${A(t) = q \cdot t - q(\xi \wedge t) = q (t - \xi)^+}$.

The expectation of $A(s)$ for $s\geq 0$ is trivially obtained from the Doob--Meyer decomposition, as follows
\begin{align}
	\label{expectationA}
	\mathbb{E}[A(s)] &= - \mathbb{E}[ \mathds{1}_{\{ \xi \le s \}} - qs ] = \mathrm{e}^{- qs} - 1 + qs.
\end{align}
Let us now compute the variance of $A(s)$ for $s\geq 0$.   On the one hand note that
\[
A(s) = q(s - \xi)^+ = q(s - \xi \wedge s), \text{ then } \operatorname{Var}[A(s)] = \operatorname{Var}[q (\xi \wedge s)].
\]
On the other hand, we have
\[
\mathbb{E}[\langle M \rangle_s] = \operatorname{Var}(M(s)) = \operatorname{Var}[ \mathds{1}_{\{\xi \le s\}}] + \operatorname{Var}[q (\xi \wedge s)] - 2 \operatorname{Cov} \Big(\mathds{1}_{\{\xi \le s\}}, q (\xi \wedge s) \Big).
\]
Hence, 
\begin{align*}
	\operatorname{Var}[A(s)] &= \mathbb{E}[\langle M \rangle_s] - \operatorname{Var}[\mathds{1}_{\{\xi \le s\}}] + 2 \operatorname{Cov} \Big(\mathds{1}_{\{\xi \le s\}}, q (\xi \wedge s) \Big). 
\end{align*}
Since   $ \mathbb{E}[\mathds{1}_{\{\xi \le s\}}] = \mathbb{E}[q (\xi \wedge s)] = 1 - \mathrm{e}^{-qs}$, we have 
\begin{align*}
	\mathbb{E}[\langle M \rangle_s]  &= 1-\mathrm{e}^{-qs},\\
	\operatorname{Var}[\mathds{1}_{\{\xi \le s\}}]& =\mathrm{e}^{-qs}(1-\mathrm{e}^{-qs}),\\
	\operatorname{Cov} \Big(\mathds{1}_{\{\xi \le s\}}, q (\xi \wedge s) \Big)&= \mathbb{E}[q (\xi \wedge s) ] - \mathbb{E}[q s \mathds{1}_{\{\xi > s\}}]- \mathbb{E}[\mathds{1}_{\{\xi \le s\}}]\mathbb{E}[q (\xi \wedge s)]\\
	&= 1-\mathrm{e}^{-qs} - qs \mathrm{e}^{-qs} - (1-\mathrm{e}^{-qs})^2 = \mathrm{e}^{-qs}(1-qs-\mathrm{e}^{-qs}).
\end{align*}
Therefore,  
\begin{equation}
	\operatorname{Var}[A(s)]    
	= 1 - \mathrm{e}^{- 2 qs} - 2 qs \mathrm{e}^{- qs}. \label{varianceA}
\end{equation}
It remains to compute the variance of $A(s_2)-A(s_1)$ for two positive reals $s_1<s_2$. \\ Using that   
$ A(s_1)(A(s_2) - A(s_1) )  = q (s_2 - s_1) A(s_1)$,
we obtain: 
\[
\mathbb{E} [ A(s_1)A(s_2) ] = q(s_2 -s_1) \mathbb{E}[ A(s_1)] + \mathbb{E} [ A(s_1)^2 ].
\]
Therefore, $\operatorname{Var}(A(s_2)-A(s_1))$ can be deduced from the expressions of  $\mathbb{E}(A(s))$  and  $\operatorname{Var}(A(s))$ obtained in \eqref{expectationA} and \eqref{varianceA} respectively. \qed

\paragraph{Acknowledgments.}
We wish to thank the anonymous reviewers for their careful reading, and several helpful remarks and suggestions which improved clarity of this manuscript.
We also thank Nathanaël Enriquez and Gabriel Faraud for fruitful discussions.

\bibliographystyle{amsalpha}
\bibliography{biblio}

\end{document}